\documentclass[aap,MSNbibl,seceqn,nameyear,dvips]{arximspdf}
\usepackage{mathrsfs}
\usepackage{graphicx}

\doi{10.1214/12-AAP846} 
\volume{23}
\issue{1}
\pubyear{2013}
\firstpage{386}
\lastpage{425}

\makeatletter

\newcommand{\dd}{\mathrm{d}}

\newtheorem{theorem}{Theorem}[section]
\newtheorem{lemma}{Lemma}[section]
\newtheorem{proposition}{Proposition}[section]
\newtheorem{corollary}{Corollary}[section]

\newcommand{\argmax}{\mathop{\arg\max}}

\makeatother

\begin{document}
\begin{frontmatter}

\title{Alpha-diversity processes and normalized inverse-Gaussian diffusions}
\runtitle{Normalized inverse-Gaussian diffusions}

\begin{aug}
\author[A]{\fnms{Matteo} \snm{Ruggiero}\corref{}\thanksref{t1}\ead[label=e1]{matteo.ruggiero@unito.it}},
\author[B]{\fnms{Stephen G.} \snm{Walker}\ead[label=e2]{S.G.Walker@kent.ac.uk}}
\and
\author[A]{\fnms{Stefano} \snm{Favaro}\thanksref{t1}\ead[label=e3]{stefano.favaro@unito.it}}
\runauthor{M. Ruggiero, S. G. Walker and S. Favaro}
\affiliation{University of Torino, University of Kent and University
of Torino}
\address[A]{M. Ruggiero\\
S. Favaro\\
Department of Economics\\
\quad and Statistics\\
University of Torino\\
Corso Unione Sovietica 218/bis\\
10134, Torino\\
Italy\\
\printead{e1}\\
\hphantom{E-mail: }\printead*{e3}}
\address[B]{S. G. Walker\\
Institute of Mathematics, Statistics\\
\quad and Actuarial Science\\
University of Kent\\
CT2 7NZ, Canterbury\\
United Kingdom\\
\printead{e2}}
\end{aug}

\thankstext{t1}{Supported by the European Research Council (ERC)
through StG ``N-BNP'' 306406.
Also affiliated with Collegio Carlo Alberto, Moncalieri, Italy.}

\received{\smonth{9} \syear{2010}}
\revised{\smonth{1} \syear{2012}}

%
\begin{abstract}
The infinitely-many-neutral-alleles model has recently been extended to
a class of diffusion processes associated with Gibbs partitions of
two-parameter Poisson--Dirichlet type. This paper introduces a family of
infinite-dimensional diffusions associated with a different subclass of
Gibbs partitions, induced by normalized inverse-Gaussian random
probability measures. Such diffusions describe the evolution of the
frequencies of infinitely-many types together with the dynamics of the
time-varying mutation rate, which is driven by an $\alpha$-diversity
diffusion. Constructed as a dynamic version, relative to this
framework, of the corresponding notion for Gibbs partitions, the latter
is explicitly derived from an underlying population model and shown to
coincide, in a special case, with the diffusion approximation of a
critical Galton--Watson branching process. The class of
infinite-dimensional processes is characterized in terms of its
infinitesimal generator on an appropriate domain, and shown to be the
limit in distribution of a certain sequence of Feller diffusions with
finitely-many types. Moreover, a discrete representation is provided by
means of appropriately transformed Moran-type particle processes, where
the particles are samples from a normalized inverse-Gaussian random
probability measure. The relationship between the limit diffusion and
the two-parameter model is also discussed.
\end{abstract}

%

\end{frontmatter}

\section{Introduction}\label{intro}

Considerable attention has been devoted recently to a class of
diffusion processes which extends the infinitely-many-neutral-alleles
model to the case of two parameters. This family takes values in the space
%
\begin{equation}\label{nabla-intro}
\overline\nabla_{\infty}=\Biggl\{z=(z_{1},z_{2},\ldots)\dvtx z_{1}\ge
z_{2}\ge\cdots\ge0, \sum_{i=1}^{\infty}z_{i}\le1\Biggr\},
\end{equation}
namely, the closure in $[0,1]^{\infty}$ of the infinite-dimensional
ordered simplex,
and is characterized, for constants $0\le\alpha<1$ and $\theta
>-\alpha
$, by the second order differential operator
%
\begin{equation}\label{operatortheta-sigma}
\mathcal{L}^{\theta,\alpha}
=\frac{1}{2}\sum_{i,j=1}^{\infty}z_{i}(\delta_{ij}-z_{j})\,\frac
{\partial
^{2}}{\partial z_{i}\,\partial z_{j}}-
\frac{1}{2}\sum_{i=1}^{\infty}(\theta z_{i}+\alpha)\,\frac{\partial
}{\partial z_{i}}
\end{equation}
acting on a certain dense sub-algebra of the space $C(\overline\nabla
_{\infty})$ of continuous functions on $\overline\nabla_{\infty}$
(throughout the paper $\delta_{ij}$ denotes Kronecker delta).
The diffusion with operator (\ref{operatortheta-sigma})
describes the evolution of the allelic frequencies at a particular
locus in a large population subject to random genetic drift and
mutation, where mutation is jointly driven by the parameters $(\theta
,\alpha)$. \citet{EK81} characterized the corresponding process when
$\alpha=0$, whereas the two-parameter family was introduced by \citet
{P09} and further investigated by \citet{RW09} and \citet{FS10}. The
latter is known to be stationary, reversible and ergodic with respect
to the Poisson--Dirichlet distribution with parameters $(\theta,\alpha
)$. This was introduced by \citet{P95} [see also \citet{P96} and \citet
{PY97}] and extends the Poisson--Dirichlet distribution of \citet{K75} as
follows. Consider a random sequence $(V_{1},V_{2},\ldots)$ obtained by
means of the so-called stick-breaking scheme
%
\begin{equation}\label{GEM}\qquad
V_{1}=W_{1},\qquad
V_{n}=W_{n}\prod_{i=1}^{n-1}(1-W_{i}),\qquad
W_{i}\stackrel{\mathrm{ind}}{\sim}\operatorname{Beta}(1-\alpha,\theta
+i\alpha),
\end{equation}
where $0\le\alpha<1$ and $\theta>-\alpha$.
The vector $(V_{1},V_{2},\ldots)$ is said to have the GEM distribution
with parameters $(\theta,\alpha)$, while the vector of descending order
statistics $(V_{(1)},V_{(2)},\ldots)$ is said to have the
Poisson--Dirichlet distribution with parameters $(\theta,\alpha)$.
The latter is also the law of the ranked frequencies of an infinite
partition induced by a two-parameter Poisson--Dirichlet random
probability measure, which generalizes the Dirichlet process introduced
by \citet{F73}. Two-parameter Poisson--Dirichlet models have found
applications in several fields. See, for example, the monographs by
\citet{B06} for fragmentation and coalescent theory, \citet{P06} for
excursion theory and combinatorics, \citet{TJ09} for machine learning,
\citet{LP09} for Bayesian inference and \citet{F10} for population
genetics. See also \citet{B08}, \citet{H09} and \citet{FLMP09}.

The Poisson--Dirichlet distribution and its two parameter extension in
turn belong to a larger class of random discrete distributions induced
by infinite partitions of Gibbs type. These were introduced by \citet
{GP05}, and applications include fragmentation and coalescent theory
[\citet{B06}, \citet{MPW08}, \citet{GMS08}], excursion theory [\citet{P03}], statistical
physics [\citet{BP07}] and Bayesian nonparametric inference [Lijoi, Mena and Pr{\"u}nster
(\citeyear{LMP05,LMP07a,LMP07b}), \citet{LPW08AAP}]. See \citet{P06} for a comprehensive
account. See also \citet{GS07}, \citet{LPW08SS} and \citet{HJL07}.

This paper introduces a class of infinite-dimensional diffusions
associated with a different subclass of Gibbs-type partitions, induced
by normalized inverse-Gaussian random probability measures. Such
discrete distributions, recently investigated by \citet{LMP05}, are
special cases of generalized gamma processes [\citet{P03}, \citet{LMP07b}], and
their intersection with two-parameter Poisson--Dirichlet models is given
by the sole case $(\theta,\alpha)=(0,1/2)$, which corresponds to a
normalized stable process with parameter $1/2$ [\citet{K75}].
The class of diffusions studied in this paper is characterized in terms
of the second order differential operator
%
\begin{eqnarray}\label{operatorgibbs}
\mathcal{A}
&=&\frac{\beta}{s}\,\frac{\partial}{\partial s}+\frac{1}{2}s\,\frac
{\partial
^{2}}{\partial s^{2}}
+\frac{1}{2}\sum_{i,j=1}^{\infty}z_{i}(\delta_{ij}-z_{j})\,\frac
{\partial
^{2}}{\partial z_{i}\,\partial z_{j}}\nonumber\\[-8pt]\\[-8pt]
&&{}-\frac{1}{2}\sum_{i=1}^{\infty}\biggl(\frac{\beta}{s}z_{i}+\alpha
\biggr)\,\frac{\partial}{\partial z_{i}}\nonumber
\end{eqnarray}
acting on a dense sub-algebra of $C_{0}([0,\infty)\times\overline
\nabla
_{\infty})$, the space of continuous functions on $[0,\infty)\times
\overline\nabla_{\infty}$ vanishing at infinity, for parameters
$(\beta
,\alpha)$, with $\beta=a\tau^{\alpha}/\alpha$, $a>0$, $\tau>0$ and
$\alpha=1/2$. By comparison with (\ref{operatortheta-sigma}), it can
be seen that the last two terms of (\ref{operatorgibbs}) describe the
time evolution of the frequencies of infinitely-many types. Common
features between (\ref{operatortheta-sigma}) and (\ref{operatorgibbs})
are the variance--covariance terms $z_{i}(\delta_{ij}-z_{j})$
and the structure of the drift or mutation terms $-[(\beta
/s)z_{i}+\alpha]$.
The distinctive
feature with respect to (\ref{operatortheta-sigma}) is
given by the fact that the positive coefficient $\theta_{t}=\beta
/S_{t}$ varies in time, and is driven by what is termed here $\alpha
$-diversity diffusion, whose operator is given by the first two terms
of (\ref{operatorgibbs}). Equivalently, $S_{t}$ follows the
stochastic differential equation
%
\begin{equation}\label{SDE}
\dd S_{t}=\frac{\beta}{S_{t}}\,\dd t+\sqrt{S_{t}}\,\dd B_{t},\qquad
S_{t}\in[0,\infty),
\end{equation}
where $B_{t}$ is a standard Brownian motion. This can be seen as a
particular instance of a continuous-time analog of the notion of
$\alpha
$-diversity, introduced by \citet{P03} for Poisson--Kingman models, which
include Gibbs-type partitions. An exchangeable random partition of
$\mathbb{N}$ is said to have $\alpha$-diversity $S$ if and only if
there exists a random variable $S$, with $0<S<\infty$ almost surely,
such that, as $n\rightarrow\infty$,
%
\begin{equation}\label{alpha-div}
K_{n}/n^{\alpha}\rightarrow S \qquad\mbox{a.s.},
\end{equation}
where $K_{n}$ is the number of classes of the partition of $\{1,\ldots
,n\}$. The connection between (\ref{SDE}) and (\ref{alpha-div}) will
become clear in Section \ref{sections-diffusion}, where the $\alpha
$-diversity diffusion will be explicitly derived.

It is to be noted that (\ref{operatortheta-sigma}) is not a special
case of (\ref{operatorgibbs}). Indeed, the only way of making
$\theta
_{t}=\beta/S_{t}$ constant is to impose null drift and volatility in
(\ref{SDE}), which implies $\theta_{t}\equiv0$. Hence, consistently
with the above recalled relation between normalized inverse-Gaussian
and Poisson--Dirichlet random measures, (\ref{operatortheta-sigma})
and (\ref{operatorgibbs}) share only the case $(\theta_{t},\alpha
)\equiv(0,1/2)$. Nonetheless, an interesting connection between these
classes of diffusions can be stated.
In particular, it will be shown that performing the same conditioning
operation in a pre-limit particle construction of normalized
inverse-Gaussian diffusions yields a particular instance of the
two-parameter model. 

The paper is organized as follows. Section \ref{sectionpreliminaries}
recalls all relevant definitions, among which are Gibbs-type
partitions, the associated generalized P\'olya-urn scheme and random
probability measures of generalized gamma and normalized
inverse-Gaussian types. Section \ref{sectionGG-result} derives some
new results on generalized gamma processes which are crucial for the
construction. These are concerned with the convergence of the number of
species represented only once in the observed sample and with the
second order approximation of the weights of the generalized P\'
olya-urn scheme associated with normalized inverse-Gaussian processes.
In Section~\ref{sections-diffusion}, by postulating simple population
dynamics underlying the time change of the species frequencies, we
derive the $\alpha$-diversity diffusion for the normalized
inverse-Gaussian case, by means of a time-varying analog of (\ref
{alpha-div}) with the limit intended in distribution, and highlight its
main properties.
In Section \ref{sectionCharacterization} normalized inverse-Gaussian
diffusions are characterized in terms of the operator (\ref
{operatorgibbs}), whose closure is shown to generate a Feller semigroup on
$C_{0}([0,\infty)\times\overline\nabla_{\infty})$, and the associated
family of processes is shown to be the limit in distribution of certain
Feller diffusions with finitely-many types.
Section \ref{sectionpopulation} provides a discrete representation of
normalized inverse-Gaussian diffusions, which are obtained as limits in
distribution of certain appropriately transformed Moran-type particle
processes which model individuals explicitly, jointly with the varying
population heterogeneity. Finally, Section \ref{secconditioning}
shows that conditioning on the $\alpha$-diversity process to be
constant, that is, $S_{t}\equiv s$, in a pre-limit version of the
particle construction yields, in the limit, the two-parameter model
(\ref{operatortheta-sigma}) with $(\theta,\alpha)=(s^{2}/4,1/2)$.

\section{Preliminaries}\label{sectionpreliminaries}

The Poisson--Dirichlet distribution and its two parameter extension
belong to the class of random discrete distributions induced by
infinite partitions of Gibbs type,
introduced by \citet{GP05}.
An exchangeable random partition of the set of natural numbers is said
to have Gibbs form
if for any $1\le k\le n$ and any $(n_{1},\ldots,n_{k})$ such that
$n_{j}\in\{1,\ldots,n\}$, for $j=1,\ldots,k$, and\vadjust{\goodbreak} $\sum
_{j=1}^{k}n_{j}=n$, the law $\Pi_{k}^{(n)}$ of the partition
$(n_{1},\ldots,n_{k})$ can be written as the product
%
\begin{equation}\label{gibbs-eppf}
\Pi_{k}^{(n)}(n_{1},\ldots,n_{k})=V_{n,k}\prod
_{j=1}^{k}(1-\alpha)_{n_{j}-1}.
\end{equation}
Here $0\le\alpha<1$,
%
\begin{equation}\label{Pochhammer}
(a)_{0}=1,\qquad
(a)_{m}=a(a+1)\cdots(a+m-1),\qquad
m>1,
\end{equation}
is the Pochhammer symbol
and the coefficients $\{V_{n,k}\dvtx k=1,\ldots,n; n\ge1\}$ satisfy the
recursive equation
%
\begin{equation}\label{V-n,k-recursion}
V_{n,k}=(n-\alpha k)V_{n+1,k}+V_{n+1,k+1}.
\end{equation}
The law of an exchangeable partition is uniquely determined by the
function $\Pi_{k}^{(n)}(n_{1},\ldots,n_{k})$, called the
\textit{exchangeable partition probability function}, which satisfies certain
consistency conditions, which imply invariance under permutations of $\{
1,\ldots,n\}$ and coherent marginalization over the $(n+1)$th item.
Hence, the law of a Gibbs partition is uniquely determined by the
family $\{V_{n,k}\dvtx k=1,\ldots,n; n\ge1\}$. Furthermore, a random
discrete probability measure governing a sequence of exchangeable
observations is said to be of a Gibbs type if it induces a partition
which can be expressed as in (\ref{gibbs-eppf}). These have associated
predictive distributions which generalize the \citet{BM73} P\'olya-urn
scheme to
%
\begin{eqnarray}\label{Pitman-urn}
&&
\mathbb{P}\{X_{n+1}\in\cdot|X_{1},\ldots,X_{n}\}\nonumber\\[-8pt]\\[-8pt]
&&\qquad=
g_{0}(n,K_{n})\nu_0(\cdot)
+g_{1}(n,K_{n})\sum
_{j=1}^{K_{n}}(n_j-\alpha
) \delta_{X_j^*}(\cdot),\nonumber
\end{eqnarray}
where $\nu_{0}$ is a nonatomic probability measure, $X_{1}^{*},\ldots
,X_{K_{n}}^{*}$ are the $K_{n}$ distinct values observed in
$X_{1},\ldots,X_{n}$ with absolute frequencies $n_{1},\ldots
,n_{K_{n}}$, and the coefficients $g_{0}$ and $g_{1}$ are given by
%
\begin{equation}\label{g0-g1}
g_{0}(n,k)=\frac{V_{n+1,k+1}}{V_{n,k}},\qquad
g_{1}(n,k)=\frac{V_{n+1,k}}{V_{n,k}}
\end{equation}
with $\{V_{n,k}\dvtx k=1,\ldots,n; n\ge1\}$ as above.
It will be of later use to note that integrating both sides of (\ref
{Pitman-urn}) yields
%
\begin{equation}\label{constraint}
g_{0}(n,K_{n})+(n-\alpha K_{n})g_{1}(n,K_{n})=1,
\end{equation}
also obtained from (\ref{V-n,k-recursion}) and (\ref{g0-g1}).
Examples of Gibbs-type random probability measures are the Dirichlet
process [\citet{F73}], obtained, for example, from (\ref{Pitman-urn}) by
setting $\theta>0$ and $\alpha=0$ in
%
\begin{equation}\label{dirichlet-coefficients}
g_{0}(n,k)=\frac{\theta+\alpha k}{\theta+n},\qquad
g_{1}(n,k)=\frac{1}{\theta+n},
\end{equation}
the two-parameter Poisson--Dirichlet process \mbox{[Pitman (\citeyear{P95,P96})]}, obtained
from (\ref{dirichlet-coefficients}) with\vadjust{\goodbreak} $0<\alpha<1$ and $\theta
>-\alpha$, the normalized stable process [\citet{K75}], obtained from
(\ref{dirichlet-coefficients}) with $0<\alpha<1$ and $\theta=0$, the
normalized inverse-Gaussian process [\citet{LMP05}] and the normalized
generalized gamma process [\citet{P03}, \citet{LMP07b}]. See also \citet{G10} for
a Gibbs-type model with finitely-many types.

The normalized generalized gamma process is a random probability
measure with representation
%
\begin{equation}\label{rpm-series}
\mu=\sum_{i=1}^{\infty}P_{i}\delta_{X_{i}},
\end{equation}
whose weights $\{P_{i},i\in\mathbb{N}\}$ are obtained by means of the
normalization
%
\begin{equation}\label{GG-normalisation}
P_{i}=\frac{J_{i}}{\sum_{k=1}^{\infty}J_{k}},
\end{equation}
where $\{J_{i},i\in\mathbb{N}\}$ are the points of a generalized gamma
process, introduced by \citet{B99}.
This is obtained from a Poisson random process on $[0,\infty)$ with
mean intensity
\[
\lambda(\dd s)=\frac{1}{\Gamma(1-\alpha)}\exp(-\tau
s)s^{-(1+\alpha)}
\,\dd s,\qquad
s\ge0,
\]
with $0<\alpha<1$ and $\tau\ge0$, so that if $N(A)$ is the number of
$J_{i}$'s which fall in $A\in\mathscr{B}([0,\infty))$, then $N(A)$ is
Poisson distributed with mean $\lambda(A)$. \citet{LMP07b} showed that a
generalized gamma random measure defined via (\ref{rpm-series}) and (\ref
{GG-normalisation}), denoted by $\operatorname{GG}(\beta,\alpha)$,
where $\beta
=a\tau^{\alpha}/\alpha$ with $a>0$ and $\tau>0$, induces a random
partition of Gibbs type with coefficients $g_{0}(n,K_{n})$ and
$g_{1}(n,K_{n})$ in (\ref{Pitman-urn}) given by
%
\begin{eqnarray}\label{g0-g1-GG}
g_{0}(n,k)&=&
\frac{ \alpha\sum_{i=0}^{n}{n\choose i}(-1)^{i}\beta
^{i/\alpha}\Gamma(k+1-i/\alpha;\beta)}
{ n \sum_{i=0}^{n-1}{n-1\choose i}(-1)^{i}\beta^{i/\alpha
}\Gamma(k-i/\alpha;\beta)},\nonumber\\[-8pt]\\[-8pt]
g_{1}(n,k)&=&
\frac{ \sum_{i=0}^{n}{n\choose i}(-1)^{i}\beta^{i/\alpha
}\Gamma(k-i/\alpha;\beta)}
{ n\sum_{i=0}^{n-1}{n-1\choose i}(-1)^{i}\beta^{i/\alpha
}\Gamma(k-i/\alpha;\beta)},
\nonumber%
\end{eqnarray}
%
where $\Gamma(c;x)$ denotes the upper incomplete gamma function
%
\begin{equation}\label{incomplete-gamma}
\Gamma(c;x)=\int_{x}^{\infty}s^{c-1}\exp(-s)\,\dd s.
\end{equation}
Special cases of a generalized gamma process with parameters $(\beta
,\alpha)$ are the Dirichlet process, obtained by letting $\tau=1$ and
$\alpha\rightarrow0$, the normalized stable process, obtained by
setting $\beta=0$, and the normalized inverse-Gaussian process,
obtained by setting $\alpha=1/2$.

We conclude the section with a brief discussion of the interpretation
of $\alpha$ in the context of species sampling with Gibbs-type
partitions. Suppose $K_{n}$ different species have been observed in the
first $n$ samples from (\ref{Pitman-urn}).
The probability that a further sample is an already observed species is
$g_{1}(n,K_{n})(n-\alpha K_{n})$, but this mass is not allocated
proportionally to the current frequencies. The
ratio of probabilities assigned to any pair of
species $(i,j)$ is
\[
r_{i,j}=\frac{n_{i}-\alpha}{n_{j}-\alpha}.
\]
When $\alpha\rightarrow0$, 
the probability of sampling species $i$ is proportional to the absolute
frequency $n_{i}$. However, since for $n_i > n_j$, $(n_{i}-\alpha)/(n_{j}-\alpha)$ is
increasing in $\alpha$, a value of $\alpha>0$ reallocates some
probability mass from type $j$ to type $i$, so that, for example, for
$n_{i}=2$ and $n_{j}=1$ we have $r_{i,j}=2, 3, 5$ for $\alpha=0, 0.5,
0.75$, respectively. Thus, $\alpha$ has a reinforcement effect on those
species that have higher frequency. See \citet{LMP07b} for a more
detailed treatment of this aspect.


\section{Some results on generalized gamma random measures}
\label{sectionGG-result}

In this section we investigate some properties of generalized gamma
random measures which will be used in the subsequent constructions. In
particular, these regard the convergence of the number of species
represented only once in the observed sample, and the second order
approximation of the weights of the generalized P\'olya-urn scheme
associated with normalized inverse-Gaussian processes.

Let $X_{1},\ldots,X_{n}$ be an $n$-sized sample drawn from a
generalized gamma process with parameters $(\beta,\alpha)$, let $K_{n}$
denote the number
of distinct species observed in the sample, and let $\mathbf
{N}_{n}:=(N_{1},\ldots,N_{K_{n}})$
denote the vector of absolute frequencies associated with each observed
species. The
probability distribution of the random variable $(K_{n},\mathbf
{N}_{n})$, for any $n\geq1$, $k=1,\ldots,n$ and frequencies
$(n_{1},\ldots,n_{k})$ such that $\sum_{i=1}^{k}n_{i}=n$, is provided
by \citet{LMP07b} and coincides with
%
\begin{eqnarray}\label{eppf}\qquad
&&
\mathbb{P}\bigl(K_{n}=k,\mathbf{N}_{n}=(n_{1},\ldots,n_{K_{n}})\bigr)\nonumber\\[-8pt]\\[-8pt]
&&\qquad =\frac{\alpha^{k-1}e^{\beta}\prod_{j=1}^{k}(1-\alpha
)_{(n_{j}-1)}}{\Gamma(n)}\sum_{s=0}^{n-1}\pmatrix{n-1\cr
s}(-1)^{s}\beta
^{s/\alpha}\Gamma\biggl(k-\frac{s}{\alpha};\beta\biggr),\nonumber
\end{eqnarray}
where $(1-\alpha)_{(n_{j}-1)}$ and $\Gamma(k-s/\alpha;\beta
)$ are as in (\ref{Pochhammer}) and (\ref{incomplete-gamma}), respectively.
Denote now by $M_{j,n}$ the number
of species represented $j$ times in the sample.
Then from equation 1.52 in \citet{P06} it follows that the distribution
of $\mathbf{M}_{n}:=(M_{1,n},\ldots,M_{n,n})$ is
given by
%
\begin{eqnarray}\label{samplingformula}
&&\mathbb{P}\bigl(\mathbf{M}_{n}=(m_{1,n},\ldots,m_{n,n})\bigr)\nonumber\\
&&\qquad= n!\frac{\alpha^{k-1}e^{\beta}}{\Gamma(n)}\prod_{j=1}^{n}\biggl(\frac
{(1-\alpha)_{(j-1)}}{j!}\biggr)^{m_{j,n}} \\
&&\qquad\quad\hspace*{59pt}{}\times\frac{1}{m_{j,n}!}\sum
_{s=0}^{n-1}\pmatrix{n-1\cr s}(-1)^{s}\beta^{s/\alpha}\Gamma
\biggl(k-\frac
{s}{\alpha};\beta\biggr)\nonumber
\end{eqnarray}
for any $n\geq1$, $k=1,\ldots,n$ and vector $(m_{1,n},\ldots
,m_{n,n})\in
\mathcal{M}_{n,k}$,
where
\[
\mathcal{M}_{n,k}=\Biggl\{(m_{1,n},\ldots,m_{n,n})\dvtx m_{i,n}\ge1,
\sum_{i=1}^{n}m_{i,n}=k, \sum_{i=1}^{n}im_{i,n}=n\Biggr\}.
\]
The following proposition identifies the speed of convergence of the
number of species represented once in the sample.
Denote by $\mathscr{C}(n,k,\alpha)$ the generalized factorial coefficient
%
\begin{equation}\label{gen-fact-coeff}
\mathscr{C}(n,k,\alpha)=\frac{1}{k!}\sum_{j=0}^{k}(-1)^{j}
\pmatrix{k\cr j}(-j\alpha)_{n},
\end{equation}
where $\mathscr{C}(0,0,\alpha)=1$ and $\mathscr{C}(n,0,\alpha)=0$.
See Charalambides [(\citeyear{C05}), Chapter~2] for a complete account.
%
\begin{proposition}
Under the normalized generalized gamma process with parameters $(\beta,
\alpha)$,
one has
%
\begin{eqnarray}\label{distribution}
&&
\mathbb{P}(M_{1,n}=m_{1,n}) \nonumber\\[-1.8pt]
&&\qquad= \frac{\alpha^{m_{1,n}-1}e^{\beta
}}{\Gamma(n)m_{1,n}!}\sum_{s=0}^{n-1}\pmatrix{n-1\cr s}(-1)^{s}\beta
^{s/\alpha}\nonumber\\[-1.8pt]
&&\hspace*{64pt}\qquad\quad{}\times\sum_{j=0}^{n-m_{1,n}}(-\alpha)^{j}\frac
{(n-m_{1,n}-j+1)_{(m_{1,n}+j)}}{j!}\\[-1.8pt]
&&\hspace*{64pt}\qquad\quad\hspace*{38pt}{}\times\sum_{k=0}^{n-m_{1,n}-j}\mathscr{C}(n-m_{1,n}-j,k,\alpha
)\nonumber\\[-1.8pt]
&&\qquad\quad\hspace*{152pt}{}\times\Gamma\biggl(k+m_{1,n}+j-\frac{s}{\alpha};\beta\biggr).\nonumber
\end{eqnarray}
Moreover,
%
\begin{equation}\label{asympdistribution}
\frac{M_{1,n}}{n^{\alpha}}\rightarrow\alpha S_{\alpha} \qquad\mbox{a.s.},
\end{equation}
where $S_{\alpha}$ is a strictly positive and almost surely finite
random variable with density function
\[
g_{S_{\alpha}}(s;\alpha,\beta)=e^{\beta-(\beta/s)^{1/\alpha
}}\frac
{f(s^{-1/\alpha};\alpha)}{\alpha s^{1+1/\alpha}}
\]
with $f(\cdot;\alpha)$ being the density of a positive stable random
variable with parameter~$\alpha$.
\end{proposition}
\begin{pf}
Denote $(x)_{[m]}=x(x-1)\cdots(x-m+1)$.
From (\ref{samplingformula}), for any $r\geq1$
one has
\begin{eqnarray*}
&&\mathbb{E}\bigl[(M_{1,n})_{[r]}\bigr]\\[-1.8pt]
&&\qquad=\sum_{k=1}^{n}\sum_{\mathcal{M}_{n,k}}n!\frac{\alpha
^{k-1}e^{\beta
}}{\Gamma(n)}\prod_{j=1}^{n}\biggl(\frac{(1-\alpha)_{(j-1)}}{j!}
\biggr)^{m_{j,n}}\frac{1}{m_{j,n}!}(m_{1,n})_{[r]}\\[-1.8pt]
&&\hspace*{16.5pt}\qquad\quad{}
\times\sum_{s=0}^{n-1}\pmatrix{n-1\cr s}(-1)^{s}\beta^{s/\alpha
}\Gamma
\biggl(k-\frac{s}{\alpha};\beta\biggr)\\
&&\qquad=\sum_{k=1}^{n}\sum_{\mathcal{M}_{n,k}}n!\frac{\alpha
^{k-1}e^{\beta
}}{\Gamma(n)(m_{1,n}-r)!} \prod_{j=2}^{n}\biggl(\frac{(1-\alpha
)_{(j-1)}}{j!}\biggr)^{m_{j,n}}\frac{1}{m_{j,n}!}\\
&&\hspace*{16.5pt}\qquad\quad{}
\times\sum_{s=0}^{n-1}\pmatrix{n-1\cr s}(-1)^{s}\beta^{s/\alpha
}\Gamma
\biggl(k-\frac{s}{\alpha};\beta\biggr)\\
&&\qquad=\sum_{k=1}^{n}\sum_{\mathcal{M}_{n-r,k-r}}n!\frac{\alpha
^{k-1}e^{\beta
}}{\Gamma(n)}\prod_{j=1}^{n}\biggl(\frac{(1-\alpha)_{(j-1)}}{j!}
\biggr)^{m_{j,n}}\frac{1}{m_{j,n}!}\\
&&\hspace*{16.5pt}\qquad\quad{}
\times\sum_{s=0}^{n-1}\pmatrix{n-1\cr s}(-1)^{s}\beta^{s/\alpha
}\Gamma
\biggl(k-\frac{s}{\alpha};\beta\biggr).
\end{eqnarray*}
In particular, by using the definition of generalized factorial
coefficient in terms of sum over the set of partitions $\mathcal
{M}_{n,k}$ [see \citet{C05}, equation 2.62], we have
\begin{eqnarray*}
&&\sum_{M_{n-r,k-r}}\prod_{j=1}^{n}\biggl(\frac{(1-\alpha
)_{(j-1)}}{j!}\biggr)^{m_{j,n}}\frac{1}{m_{j,n}!}\\
&&\qquad=\frac
{(n)_{[r]}}{n!\alpha^{k-r}}\mathscr{C}(n-r,k-r,\alpha).
\end{eqnarray*}
Therefore, we
obtain
%
\begin{eqnarray}\label{factorialmoment}
&&\mathbb{E}\bigl[(M_{1,n})_{[r]}\bigr]\nonumber\\
&&\qquad=\sum_{k=1}^{n}\frac{\alpha
^{r-1}(n)_{[r]}e^{\beta}}{\Gamma(n)
}\mathscr{C}(n-r,k-r,\alpha)\\
&&\qquad\quad\hspace*{14.5pt}{}
\times\sum_{s=0}^{n-1}\pmatrix{n-1\cr s}(-1)^{s}\beta^{s/\alpha
}\Gamma
\biggl(k-\frac{s}{\alpha};\beta\biggr).\nonumber
\end{eqnarray}
In order to obtain the distribution of the random variable $M_{1,n}$,
we can make use of the probability generating function of $M_{1,n}$,
denoted $G_{(M_{1,n})}(t)$. From (\ref{factorialmoment}) we have
\begin{eqnarray*}
G_{(M_{1,n})}(t)
&=&\sum_{r=0}^{\infty}\frac{\alpha^{r-1}e^{\beta}(n)_{[r]}}{\Gamma
(n)}\\[-2pt]
&&\hspace*{13.3pt}{} \times\sum_{s=0}^{n-1}\pmatrix{n-1\cr s}(-1)^{s}\beta^{s/\alpha}\\[-2pt]
&&\hspace*{41pt}{} \times\sum_{k=0}^{n}\mathscr
{C}(n-r,k,\alpha)\Gamma\biggl(k+r-\frac{s}{\alpha};\beta\biggr)\frac
{(t-1)^{r}}{r!}.
\end{eqnarray*}
Therefore, the distribution of $M_{1,n}$ is given by
\begin{eqnarray*}
&&
\mathbb{P}(M_{1,n}=m_{1,n})\\[-2pt]
&&\qquad=\frac{1}{m_{1,n}!}\sum_{j=0}^{\infty
}\frac
{\alpha^{m_{1,n}+j-1}e^{\beta}(n)_{[m_{1,n}+j]}}{\Gamma(n)}\\[-2pt]
&&\hspace*{39pt}\qquad\quad{}\times\sum
_{s=0}^{n-1}\pmatrix{n-1\cr s}(-1)^{s}\beta^{s/\alpha}\\[-2pt]
&&\hspace*{64pt}\qquad\quad{}\times\sum_{k=0}^{n}\mathscr{C}(n-m_{1,n}-j,k,\alpha)\\[-2pt]
&&\hspace*{87pt}\qquad\quad{}\times\Gamma
\biggl(k+m_{1,n}+j-\frac{s}{\alpha};\beta\biggr)\frac{d^{m_{1,n}}}{dt^{m_{1,n}}}
\frac{(t-1)^{m_{1,n}+j}}{(m_{1,n}+j)!}\bigg|_{t=0}\\[-2pt]
&&\qquad=\frac{1}{m_{1,n}!}\sum_{j=0}^{\infty}\frac{\alpha
^{m_{1,n}+j-1}e^{\beta}(n)_{[m_{1,n}+j]}}{\Gamma(n)}\\[-2pt]
&&\hspace*{39pt}\qquad\quad{}\times\sum
_{s=0}^{n-1}\pmatrix{n-1\cr s}(-1)^{s}\beta^{s/\alpha}\\[-2pt]
&&\hspace*{64pt}\qquad\quad{}\times\sum_{k=1}^{n-m_{1,n}-j}\mathscr{C}(n-m_{1,n}-j,k,\alpha)\\[-2pt]
&&\hspace*{113pt}\qquad\quad{}\times\Gamma
\biggl(k+m_{1,n}+j-\frac{s}{\alpha};\beta\biggr)\frac{(-1)^{j}}{j!}\\[-2pt]
&&\qquad=\frac{\alpha^{m_{1,n}-1}e^{\beta}}{\Gamma(n)m_{1,n}!}\sum
_{s=0}^{n-1}\pmatrix{n-1\cr s}(-1)^{s}\beta^{s/\alpha}\\[-2pt]
&&\hspace*{60.6pt}\qquad\quad{}\times\sum_{j=0}^{n-m_{1,n}}(-\alpha)^{j}\frac
{(n-m_{1,n}-j+1)_{(m_{1,n}+j)}}{j!}\\[-2pt]
&&\hspace*{98pt}\qquad\quad{}\times\sum_{k=0}^{n-m_{1,n}-j}\mathscr{C}(n-m_{1,n}-j,k,\alpha)\\[-2pt]
&&\hspace*{148pt}\qquad\quad{}\times\Gamma
\biggl(k+m_{1,n}+j-\frac{s}{\alpha};\beta\biggr),
\end{eqnarray*}
where the last identity is due to the fact that $\mathscr
{C}(n,k,\alpha
)=0$ for any $k>n$. Proposition 3 in \citet{LMP07b} shows that
%
\begin{equation}\label{GG-alpha-diversity}
K_{n}/n^{\alpha}\rightarrow S_{\alpha}
\end{equation}
almost surely, where
$S_{\alpha}$ is an almost surely positive and finite random variable
with density function
\[
g_{S_{\alpha}}(s;\alpha,\beta)=e^{\beta-(\beta/s)^{1/\alpha
}}\frac
{f(s^{-1/\alpha};\alpha)}{\alpha s^{1+1/\alpha}}
\]
with $f(\cdot;\alpha)$ being the density function of a positive stable
random variable with parameter $\alpha$. In other terms,
according to Definition 3.10 in \citet{P06}, an exchangeable partition
of $\mathbb{N}$
having EPPF (\ref{eppf}) has $\alpha$-diversity $S_{\alpha}$. A simple
application of Lemma 3.11
in \citet{P06} leads to (\ref{asympdistribution}).
\end{pf}

A second aspect of generalized gamma random measures we need to address
for later use is the approximate behavior of the coefficients in the
generalized P\'olya urn (\ref{Pitman-urn}).
It is well known that the first order behavior of (\ref{g0-g1-GG}) is
that of a normalized stable process, that is,
%
\begin{equation}\label{pesi-stable-case}
g_{0}(n,k)\approx\alpha k/n,\qquad
g_{1}(n,k)\approx1/n,
\end{equation}
also implied by the next result.
However, it turns out that for the definition of the diffusion
processes which are the object of the next two sections, it is crucial
to know the second order approximation. The following proposition,
whose proof is deferred to the \hyperref[app]{Appendix}, identifies such behavior for
the normalized inverse-Gaussian case $\alpha=1/2$.
%
\begin{proposition}\label{second-order-approximation}
Let $g_{0}(n,k)$ and $g_{1}(n,k)$ be as in (\ref{g0-g1-GG}). When
$\alpha=1/2$,
\[
g_{0}(n,k)=\frac{\alpha k}{n}+\frac{\beta/s_{n}}{n}+o(n^{-1})
\]
and
%
\begin{equation}\label{g1-approximation2}
g_{1}(n,k)=\frac{1}{n}-\frac{\beta/s_{n}}{n^{2}}+o(n^{-2}),
\end{equation}
where $s_{n}=k/n^{\alpha}$ and $\beta=a\tau^{\alpha}/\alpha$.
\end{proposition}



\section{Alpha-diversity processes}\label{sections-diffusion}

Making use of the results of the previous section, here we construct a
one-dimensional diffusion process which can be seen as a dynamic
version of the notion of $\alpha$-diversity, recalled in (\ref
{alpha-div}), relative to the case of\vadjust{\goodbreak} normalized inverse-Gaussian
random probability measures. Such diffusion, which will be crucial for
the construction of Section~\ref{sectionCharacterization}, is obtained
as weak limit of an appropriately rescaled random walk on the integers,
whose dynamics are driven by an underlying population process. This is
briefly outlined here and will be formalized in Section
\ref{sectionpopulation}. Consider $n$ particles, denoted
$x^{(n)}=(x_{1},\ldots,x_{n})$ with $x_{i}\in\mathbb{X}$ for each $i$,
where $\mathbb{X}$ is a Polish space, and denote by
$K_{n}=K_{n}(x^{(n)})$ the number of distinct values observed in
$(x_{1},\ldots,x_{n})$. Let the vector $(x_{1},\ldots,x_{n})$ be
updated at discrete times by replacing a uniformly chosen coordinate.
Conditionally on $K_{n}(x^{(n)})=k$, the incoming particle will be a
copy of one still in the vector, after the removal, with probability
$g_{1}(n-1,k_{r})$, and will be a new value with probability
$g_{0}(n-1,k_{r})$, where $g_{1}(n-1,k)$ and $g_{0}(n-1,k)$ are as in
(\ref{g0-g1-GG}) and $k_{r}$ is the value of $k$ after the removal.
Denote by $\{ K_{n}(m),m\in\mathbb{N}_{0}\}$ the chain which keeps
track of the number of distinct types in $(x_{1},\ldots,x_{n})$. Then,
letting $m_{1,n}$ be the number of clusters of size one in
$(x_{1},\ldots,x_{n})$, which, by means of (\ref{asympdistribution})
and (\ref{GG-alpha-diversity}) is approximately $\alpha k$ for large
$n$, the transition probabilities for $K_{n}(m)$,
\[
p(k,k')=\mathbb{P}\{K_{n}(m+1)=k'|K_{n}(m)=k\}
\]
%
are asymptotically equivalent to
%
\begin{equation}\label{K-transition}
p(k,k')=\cases{\displaystyle
\biggl(1-\frac{\alpha k}{n}\biggr)g_{0}(n-1,k), \hspace*{42pt}\quad \mbox{if
$1 \leq k < n,
k'=k+1$},\vspace*{2pt}\cr
\displaystyle
\frac{\alpha k}{n}g_{1}(n-1,k-1)\bigl(n-1-\alpha(k-1)\bigr),\cr
\qquad\hspace*{136pt}
\mbox{if $1 < k \leq n, k'=k-1$},\vspace*{2pt}\cr
\displaystyle
1-p(k,k+1)-p(k,k-1), \qquad \mbox{if $k'=k$},\vspace*{2pt}\cr
0, \qquad\hspace*{126.2pt} \mbox{else}}\hspace*{-35pt}
\end{equation}
for $1\le k\le n$. That is, with probability $m_{1,n}/n\approx \alpha
k/n$ a cluster of size one is selected and removed, with probability
$g_{0}(n-1,k)$ a new species appears and with probability
$g_{1}(n-1,k)(n-1-\alpha(k-1))$ a survivor has an offspring. Note that
$k=1$ and $k=n$
are set to be barriers, to render the fact that $m_{1,n}$ equals 0 and~$n$ when $k$
equals 1 and $n$, respectively.


The following theorem finds the conditions under which the rescaled
chain $K_{n}(m)/n^{\alpha}$ converges to a diffusion process on
$[0,\infty)$. Here we provide a sketch of the proof with the aim of
favoring the intuition. The formalization of the result is contained in
the proof of Theorem \ref{theoremparticle-convergence}, while that of
the fact that the limiting diffusion is
well defined, that is, the corresponding operator generates a Feller
semigroup on an appropriate subspace of $C([0,\infty))$, is provided in
Corollary \ref{S-diff-well-defined} below.

Throughout the paper $C_{B}(A)$ denotes the space of continuous
functions from $A$ to $B$, while $X_{n}\Rightarrow X$ denotes
convergence in distribution.
%
\begin{theorem}\label{S-convergence}
Let $\{K_{n}(m),m\in\mathbb{N}_{0}\}$ be a Markov chain with transition
probabilities as in (\ref{K-transition})\vadjust{\goodbreak} determined by a generalized
gamma process with $\beta\ge0$ and $\alpha=1/2$, and define $\{
\tilde K_{n}(t),t\ge0\}$ to be such that $\tilde K_{n}(t)=K_{n}(\lfloor
n^{3/2}t\rfloor)/n^{\alpha}$. Let also $\{S_{t},t\ge0\}$ be a
diffusion process driven by the stochastic differential equation
%
\begin{equation}\label{SDE2}
\dd S_{t}=\frac{\beta}{S_{t}}\,\dd t+\sqrt{S_{t}} \,\dd B_{t},\qquad
S_{t}\ge0,
\end{equation}
where $B_{t}$ is a standard Brownian motion. If $\tilde K_{n}(0)\Rightarrow S_{0}$, then
%
\begin{equation}\label{weak-conv-S}
\{\tilde K_{n}(t),t\ge0\}\Rightarrow\{S_{t},t\ge0\}
\qquad\mbox{in } C_{[0,\infty)}([0,\infty))\mbox{ as }n\rightarrow
\infty.
\end{equation}
\end{theorem}
\begin{pf}
Let $\alpha=1/2$.
From Proposition \ref
{second-order-approximation} we can write (\ref{K-transition}) as
follows (for ease of presentation we use $n$ and $k$ in place of $n-1$
and $k-1$ since it is asymptotically equivalent):
\[
p(k,k')=\cases{
\displaystyle\biggl(1-\frac{\alpha k}{n}\biggr)\biggl(\frac{\alpha k}{n}+\frac{\beta
/s_{n}}{n}\biggr)+o(n^{-1}),
&\hspace*{-3.5pt}\quad if $1 \leq k < n, k'=k+1$,\vspace*{2pt}\cr
\displaystyle\frac{\alpha k}{n}\biggl(\frac{1}{n}-\frac{\beta/s_{n}}{n^{2}}
\biggr)(n-\alpha k)+o(n^{-3/2})
,&\hspace*{-3.5pt}\quad if $1 < k \leq n, k'=k-1$,\vspace*{2pt}\cr
\displaystyle1-p(k,k+1)-p(k,k-1)+o(n^{-1}),
&\hspace*{-3.5pt}\quad if $k'=k$,\cr
0,&\hspace*{-3.5pt}\quad else.}
\]
%
The conditional expected increment of the process $\{K_{n}(m)/n^{\alpha
},m\in\mathbb{N}_{0}\}$ is
%
\begin{eqnarray}\label{exp-increm}
&&\mathbb{E}\biggl(\frac{k'}{n^{\alpha}}-\frac{k}{n^{\alpha}}\Big|k\biggr)\nonumber\\
&&\qquad= \frac{1}{n^{\alpha}}\biggl[\biggl(1-\frac{\alpha k}{n}\biggr)
\biggl(\frac{\alpha k}{n}+\frac{\beta/s_{n}}{n}\biggr)
-\frac{\alpha k}{n}\biggl(\frac{1}{n}-\frac{\beta/s_{n}}{n^{2}}\biggr)
(n-\alpha
k)\biggr]\nonumber\\[-8pt]\\[-8pt]
&&\qquad\quad{}+o\biggl(\frac{1}{n^{1+\alpha}}\biggr)\nonumber\\
&&\qquad= \frac{\beta/s_{n}}{n^{1+\alpha}}
+o\biggl(\frac{1}{n^{1+\alpha}}
\biggr).\nonumber
\end{eqnarray}
Similarly, the conditional second moment of the increment is
%
\begin{eqnarray}\label{exp-variance}
&&\mathbb{E}\biggl[\biggl(\frac{k'}{n^{\alpha}}-\frac{k}{n^{\alpha}}
\biggr)^{2}\Big|k\biggr]\nonumber\\
&&\qquad= \frac{1}{n^{2\alpha}}\biggl[\biggl(1-\frac{\alpha k}{n}\biggr)
\biggl(\frac{\alpha k}{n}+\frac{\beta/s_{n}}{n}\biggr)
+\frac{\alpha k}{n}\biggl(\frac{1}{n}-\frac{\beta/s_{n}}{n^{2}}\biggr)
(n-\alpha k)\biggr]\nonumber\\[-8pt]\\[-8pt]
&&\qquad\quad{}+o\biggl(\frac{1}{n^{1+2\alpha}}\biggr)\nonumber\\
&&\qquad= \frac{2\alpha k}{n^{1+2\alpha}}+o\biggl(\frac{1}{n^{1+2\alpha}}
\biggr).\nonumber
\end{eqnarray}
Since $k\approx s n^{\alpha}$, and recalling that $s_{n}\rightarrow s$
almost surely, we have
\[
n^{1+\alpha} \mathbb{E}\biggl(\frac{k'}{n^{\alpha}}-\frac{k}{n^{\alpha
}}\Big|k\biggr)\rightarrow\beta/s
\]
and
\[
n^{1+\alpha} \mathbb{E}\biggl[\biggl(\frac{k'}{n^{\alpha}}-\frac
{k}{n^{\alpha}}\biggr)^{2}\Big|k\biggr]\rightarrow2\alpha s.
\]
It is easy to check that all conditional $m$th moments of $\Delta
k/n^{\alpha}$ converge to zero for $m\ge3$, whence it follows by
standard theory [cf., e.g., \citet{KT81}] that, as $n\rightarrow\infty$,
the process $\tilde K_{n}(t)=K_{n}(\lfloor n^{3/2}t\rfloor)/n^{\alpha}$
converges in distribution to a diffusion process $S_{t}$ on $[0,\infty
)$ with drift $\beta/S_{t}$ and diffusion coefficient $\sqrt{2\alpha S_{t}}$.
\end{pf}

As anticipated, the second order approximation of $g_{0}(n,k)$ is
crucial for establishing the drift of the limiting diffusion, as the
first order terms cancel.
It is interesting to note that when $\beta=0$, which yields the
normalized stable case, the limiting diffusion reduces to the diffusion
approximation of a critical Galton--Watson branching process, also known
as the zero-drift Feller diffusion. See, for example, \citet{EK86},
Theorem 9.1.3. This also holds approximately for high values of
$S_{t}$, in which case the drift becomes negligible.

\begin{figure}
\begin{tabular}{@{}c@{}}

\includegraphics{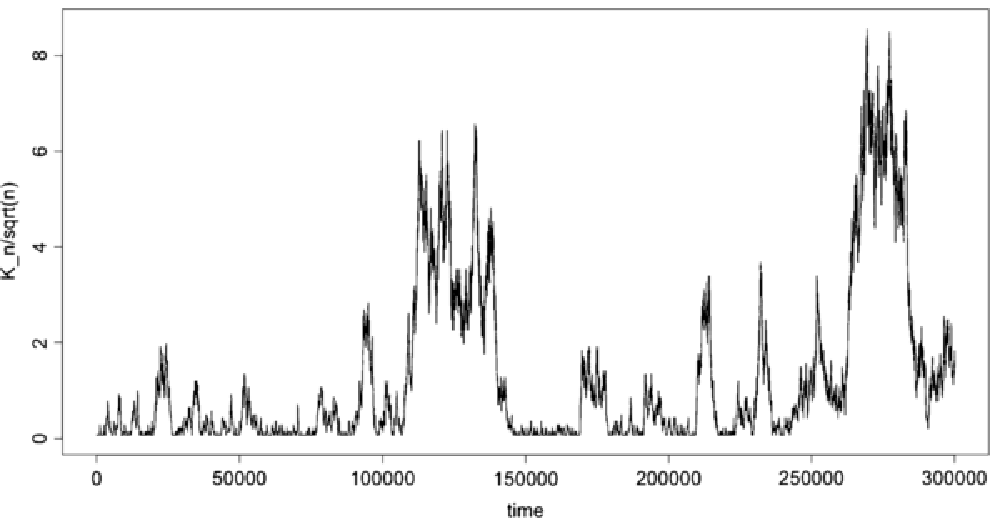}
\\
(a)\\[4pt]

\includegraphics{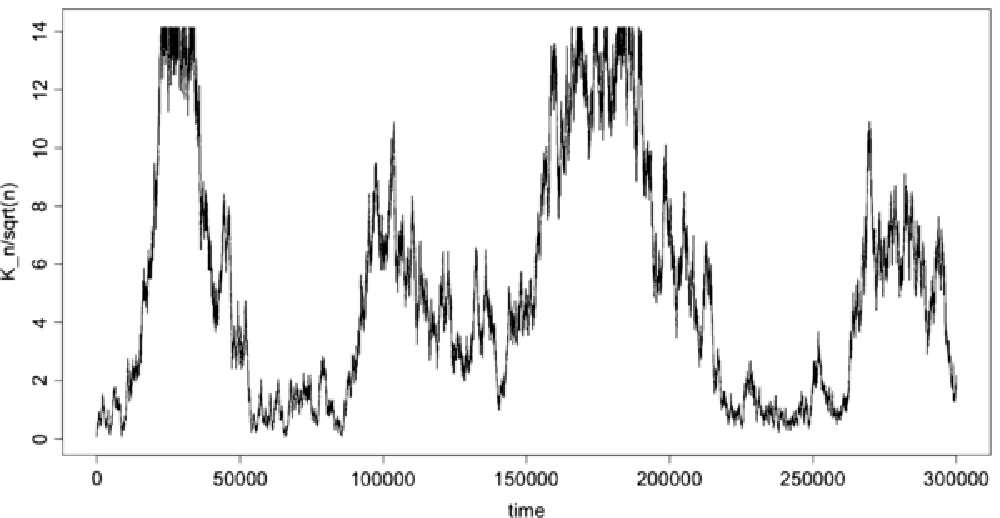}
\\
(b)\\[4pt]

\includegraphics{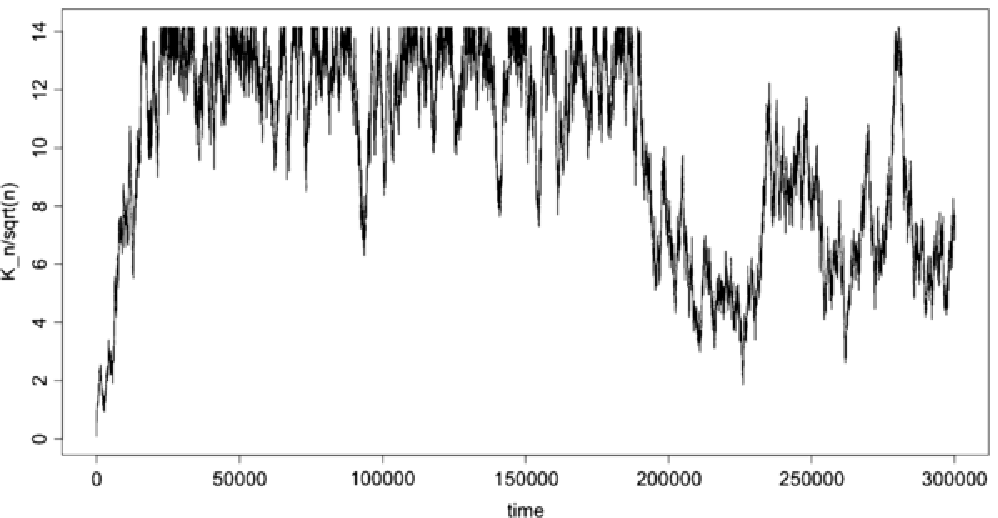}
\\
(c)
\end{tabular}
\caption{Three sample paths of the random walk $\{K_{n}(m)/n^{\alpha
},m\in\mathbb{N}_{0}\}$, with dynamics as in Theorem \protect\ref
{S-convergence}, starting from $1/\sqrt{n}$ with $n=200$, for parameter
values: \textup{(a)} $\beta=0$, \textup{(b)} $\beta=100$,
\textup{(c)}~$\beta =1000$. The figures show how $\beta$ influences the
dynamic clustering structure in the population.}\label{figpaths}
\end{figure}

In order to have some heuristics on the behavior of the $\alpha
$-diversity process, Figure \ref{figpaths} shows $3\times10^{5}$
steps of the random walk $\{K_{n}(m)/n^{\alpha},m\in\mathbb{N}_{0}\}$
with dynamics as in Theorem \ref{S-convergence}, starting from
$1/\sqrt
{n}$ with $n=200$. The three paths correspond to $\beta$ being equal to
0, 100 and 1000. It is apparent how $\beta$ influences the dynamic
clustering structure in the population.

It is well known that when $\beta=0$, the point 0 is an absorbing
boundary for $S_{t}$.
The next result provides the boundary classification, using Feller's
terminology, for the case $\beta>0$.
%
\begin{proposition}\label{boundaries}
Let $S_{t}$ be as in Theorem \ref{S-convergence} with $\beta>0$.
Then the points 0 and $\infty$ are, respectively, an entrance and a
natural boundary.
\end{proposition}
\begin{pf}
The scale function for the process, defined as
%
\begin{equation}\label{scale-function}
S(x)=\int_{x_{0}}^{x}s(y)\,\dd y,\qquad 0<x<\infty,
\end{equation}
where
\[
s(y)=\exp\biggl\{-\int_{y_{0}}^{y}\frac{2\mu(t)}{\sigma^{2}(t)}\,\dd t
\biggr\}
\]
and $\mu(x)$ and $\sigma^{2}(x)$ denote drift and diffusion,
equals
\begin{eqnarray*}
S(x)
&=& \int_{x_{0}}^{x}\exp\biggl\{-2\beta\biggl(\frac{1}{y_{0}}-\frac
{1}{y}\biggr)\biggr\}\,\dd y\\
&=& e^{-2\beta/y_{0}}[xe^{2\beta/x}-x_{0}e^{2\beta/x_{0}}-2\beta
\operatorname{Ei}(2\beta/x)+\operatorname{Ei}(2\beta/x_{0})],
\end{eqnarray*}
where $\operatorname{Ei}(z)$ is the exponential integral
\[
\operatorname{Ei}(z)=-\int_{-z}^{\infty}t^{-1}e^{-t}\,\dd t.
\]
Letting $S[a,b]=S(b)-S(a)$, for $0<a<b<\infty$, we have
%
\begin{eqnarray}\label{scale}
S(0,b]&=&\lim_{a\downarrow0}S[a,b]=\infty,\nonumber\\[-8pt]\\[-8pt]
S[a,\infty)&=&\lim_{b\uparrow\infty}S[a,b]=\infty.
\nonumber
\end{eqnarray}
Moreover, the speed measure is given by
\begin{eqnarray*}
M[c,d]
&=& \int_{c}^{d}[\sigma^{2}(t)s(t)]^{-1}\,\dd t\\
&=& e^{2\beta/y_{0}}[\operatorname{Ei}(-2\beta/c)-\operatorname
{Ei}(-2\beta
/d)]
\end{eqnarray*}
from which $M(0,d]=\lim_{c\downarrow0}M[c,d]<\infty$ and
%
\begin{equation}\label{speed}
M[c,\infty)
=\lim_{d\uparrow\infty}M[c,d]=\infty.
\end{equation}
Now (\ref{scale}) implies that
\begin{eqnarray*}
\Sigma(0) &=& \lim_{l\downarrow0}\int_{l}^{x}S(l,y]\,\dd M(y)=\infty,\\
\Sigma(\infty) &=& \lim_{l\uparrow\infty}\int_{x}^{r}S[y,r)\,\dd
M(y)=\infty
\end{eqnarray*}
and (\ref{speed}) implies
\[
N(\infty)=\lim_{r\uparrow\infty}\int_{x}^{r}S[x,y]\,\dd M(y)=\infty,
\]
while
\begin{eqnarray*}
N(0) &=& \lim_{l\downarrow0}\int_{l}^{x}S[y,x]\,\dd M(y)\\
&=& \lim_{l\downarrow0}\int_{l}^{x}\frac{e^{-2\beta/y}}{y}
\biggl[xe^{2\beta/x} - y e^{2\beta/y}+2\beta\biggl(\operatorname{Ei}\biggl(\frac
{2\beta}{y}\biggr) - \operatorname{Ei}\biggl(\frac{2\beta}{x}\biggr)\biggr)
\biggr]\,\dd y<\infty
\end{eqnarray*}
since
\[
\lim_{y\downarrow0}\frac{e^{-2\beta/y}}{y}\operatorname{Ei}\biggl(\frac
{2\beta
}{y}\biggr)<\infty.
\]
%
The statement now follows from, for example, \citet{KT81},
Section~15.6.
\end{pf}

Hence, when $\beta>0$ neither boundary point is attainable from the
interior of the state space, from which the actual state space is
$[0,\infty)$ for $\{S_{t},t\ge0\}$ and $(0,\infty)$ for $\{
S_{t},t>0\}
$. The process can be made to start at 0, in which case it instantly
moves toward the interior of the state space and never comes back.
Consequently, we will use $(0,\infty)$ or $[0,\infty)$ as state space
at convenience, with the agreement that $(0,\infty)$ is referred to $\{
S_{t},t>0\}$.

As a corollary, we formalize the well-definedness of the $\alpha
$-diversity diffusion. Denote by $C_{0}(K)$ the space of continuous
functions vanishing at infinity on a locally compact set $K$, and let
$\|\cdot\|$ be a norm which makes $C_{0}(K)$ a Banach space. Recall
that a Feller semigroup on $C_{0}(K)$ is a one-parameter family of
bounded linear operators $\{T(t), t\ge0\}$ on $C_{0}(K)$ such that
$T(t)$ has the semigroup property $T(s+t)=T(s)T(t)$ for all $s,t\ge0$,
is strongly continuous, that is,
\[
\|T(t)f-f\|\rightarrow0 \qquad\mbox{as }t\rightarrow0, f\in C_{0}(K),
\]
and, for all $t\ge0$, $T(t)$ is a contraction, that is, $\|T(t)\|\le
1$, is conservative in the sense that $T(t)1=1$, and is positive in the
sense that it preserves the cone of nonnegative functions.
%
\begin{corollary}\label{S-diff-well-defined}
For $\beta\ge0$, let $\mathcal{A}_{0}$ be the second order
differential operator
%
\begin{equation}\label{S-generator}
\mathcal{A}_{0}=\frac{\beta}{s}\,\frac{\dd}{\dd s}+\frac
{1}{2}s\,\frac{\dd
^{2}}{\dd s^{2}}
\end{equation}
and define
%
\begin{equation}\label{domain-A0}
\mathscr{D}(\mathcal{A}_{0})=\{f\in C_{0}([0,\infty))\cap
C^{2}((0,\infty))\dvtx \mathcal{A}_{0}f\in C_{0}([0,\infty))\}.
\end{equation}
Then $\{(f,\mathcal{A}_{0}f)\dvtx f\in\mathscr{D}(\mathcal{A}_{0})\}$
generates a Feller semigroup on $C_{0}([0,\infty))$.
\end{corollary}
\begin{pf}
The result follows from Proposition \ref{boundaries} together with
Corollary~8.1.2 in \citet{EK86}.
\end{pf}

%

An immediate question that arises is whether the $\alpha$-diversity
diffusion is stationary.
The following proposition, which concludes the section, provides a
negative answer.
%
\begin{proposition}
Let $\{S_{t},t\ge0\}$ be as in Theorem \ref{S-convergence}. Then there
exists no stationary density for the process.\vadjust{\goodbreak}
\end{proposition}
\begin{pf}
A stationary density, if it exists, is given by
\[
\psi(x)=m(x)[C_{1}S(x)+C_{2}],\qquad x\ge0,
\]
where $m(x)=[s(x)\sigma^{2}(x)]^{-1}$, $s(x)$ and $S(x)$ are as in
(\ref
{scale-function}), and $C_{1},C_{2}$ are constants determined in order
to guarantee the nonnegativity and integrability to one of $\psi$.
Here $s(x)=e^{-2\beta/x}$ and
\[
S(x)=xe^{2\beta/x}-2\beta\operatorname{Ei}(2\beta/x)
\]
so that
\[
\psi(x)=C_{1}-2\beta C_{1}x^{-1}e^{-2\beta/x}\operatorname
{Ei}(2\beta
/x)+C_{2}x^{-1}e^{-2\beta/x}.
\]
The second term is not integrable in a neighborhood of infinity, since
there exists an $x_{0}>0$ such that
\[
-x^{-1}e^{-2\beta/x}\operatorname{Ei}(2\beta/x)>x^{-1}
\qquad\mbox{for all }x>x_{0},
\]
hence $C_{1}$ must be zero.
Since neither the third term is integrable, this gives
the result.
\end{pf}


\section{Normalized inverse-Gaussian diffusions}
\label{sectionCharacterization}

The $\alpha$-diver\-sity process constructed in the previous section is
a key component in the definition of the class of normalized
inverse-Gaussian diffusions. In this section we characterize such
infinite-dimensional processes in terms of their infinitesimal
generator, and show that they can be obtained as the limit in
distribution of a certain sequence of Feller diffusions with
finitely-many types.
The association of the limit family with the class of normalized
inverse-Gaussian random probability measures will instead be shown in
Section \ref{sectionpopulation}.

Consider the $(n-1)$-dimensional simplex
\[
\Delta_{n}=\Biggl\{z\in[0,1]^{n}\dvtx z_{i}\ge0, \sum
_{i=1}^{n}z_{i}=1\Biggr\}
\]
and the closed subspace of $\Delta_{n}$ given by
\[
\tilde\Delta_{n}=\Biggl\{z\in[0,1]^{n}\dvtx z_{i}\ge\varepsilon_{n}, \sum
_{i=1}^{n}z_{i}=1\Biggr\},
\]
so that $\varepsilon_{n}\le z_{i}\le1-(n-1)\varepsilon_{n}$ for
$z_{i}\in\tilde\Delta_{n}$,
where $\{\varepsilon_{n}\}\subset\mathbb{R}_{+}$ is a nonincreasing
sequence such that
%
\begin{equation}\label{varepsilon-n}
0<\varepsilon_{n}<\frac{1}{n} \qquad\forall n\ge2,\qquad
n\varepsilon_{n}\downarrow0.
\end{equation}
%
Define, for $(z_{0},z_{1},\ldots,z_{n})\in(0,\infty)\times\tilde
\Delta
_{n}$, the differential operator
\[
\mathcal{A}_{n}
=\frac{1}{2}\sum_{i,j=0}^{n}a_{ij}^{(n)}(z)\,\frac{\partial
^{2}}{\partial z_{i}\,\partial z_{j}}
+\frac{1}{2}\sum_{i=0}^{n}b_{i}^{(n)}(z)\,\frac{\partial}{\partial z_{i}},
\]
where the covariance components $(a_{ij}^{(n)}(z))_{i,j=0,\ldots,n}$
are set to be
\[
a_{ij}^{(n)}(z)=
\cases{
\displaystyle z_{0},&\quad$ i=j=0 $,\vspace*{2pt}\cr
\displaystyle
(z_{i}-\varepsilon_{n})\bigl(\delta_{ij}(1-n\varepsilon
_{n})-(z_{j}-\varepsilon_{n})\bigr)
, &\quad $ 1\le i,j\le n$,\vspace*{2pt}\cr
\displaystyle
0,&\quad else,}
\]
and, for $\beta>0$, the drift components are
\begin{eqnarray*}
b_{0}^{(n)}(z) &=& \frac{\beta}{z_{0}},\\
b_{i}^{(n)}(z) &=& \frac{\beta}{z_{0}(n-1)}(1-z_{i})-\frac{\beta
}{z_{0}}z_{i}\\
&&{}-\alpha\bigl(1-\exp\{-(z_{i}-\varepsilon
_{n})e^{1/\varepsilon_{n}}\}\bigr),\qquad
i=1,\ldots,n.
\end{eqnarray*}
Observe that $a_{ij}^{(n)}(z)$, for $1\le i,j\le n$, can be seen as a
Wright--Fisher type covariance restricted to $[\varepsilon
_{n},1-(n-1)\varepsilon_{n}]^{n}$, since for such indices $i,j$
%
\begin{equation}\label{a-ij}
a_{ij}^{(n)}(z)=
\cases{\displaystyle
(z_{i}-\varepsilon_{n})\bigl(1-(n-1)\varepsilon_{n}-z_{i}\bigr),&\quad $ i=j
$,\vspace*{2pt}\cr
\displaystyle
-(z_{i}-\varepsilon_{n})(z_{j}-\varepsilon_{n}),&\quad $ i\ne j$,}
\end{equation}
and that the first two terms in $b_{i}^{(n)}(z)$, $i=1,\ldots,n$, equal
\[
\frac{\beta}{z_{0}(n-1)}\bigl(1-(n-1)\varepsilon_{n}-z_{i}\bigr)-\frac{\beta
}{z_{0}}(z_{i}-\varepsilon_{n})
\]
from which the behavior at the boundary is clear.
For ease of exposition and in analogy with the previous section, whenever
convenient we will
denote $z_{0}$ by $s$, so that for instance $\mathcal{A}_{n}$ can be
written more explicitly,
%
\begin{equation}\label{generatordelta-n}\quad
\mathcal{A}_{n}
=
\frac{\beta}{s}\,\frac{\partial}{\partial s}+\frac{1}{2}s\,\frac
{\partial
^{2}}{\partial s^{2}}
+\frac{1}{2}\sum_{i,j=1}^{n}a_{ij}^{(n)}(z) \,\frac{\partial
^{2}}{\partial z_{i}\,\partial z_{j}}+\frac{1}{2}\sum
_{i=1}^{n}b_{i}^{(n)}(z)\,\frac{\partial}{\partial z_{i}}
\end{equation}
with $a_{ij}^{(n)}(z)$ and $b_{i}^{(n)}(z)$ as above.
The domain of $\mathcal{A}_{n}$ is taken to be
%
\begin{equation}\label{domain-delta-n}
\mathscr{D}(\mathcal{A}_{n})
=\{f\dvtx f=f_{0}\times f_{1}, f_{0}\in\mathscr{D}(\mathcal{A}_{0}),
f_{1}\in C^{2}(\tilde\Delta_{n})
\},
\end{equation}
where $(f_{0}\times f_{1})(s,z)=f_{0}(s)f_{1}(z)$, $\mathscr
{D}(\mathcal
{A}_{0})$ is (\ref{domain-A0}), and
\[
C^{2}(\tilde\Delta_{n})=
\{f\in C(\tilde\Delta_{n})\dvtx \exists\tilde f\in C^{2}(\mathbb
{R}^{n}), \tilde f|_{\tilde\Delta_{n}}=f
\}.
\]
%
The operator $\mathcal{A}_{n}$ drives $n+1$ components: those labeled
from 1 to $n$ can be seen as the frequencies associated to $n$ species
in a large population, bounded from below by $\varepsilon_{n}$; the
$z_{0}$ or $s$ component is a positive real variable which evolves
independently according to the $\alpha$-diversity diffusion (\ref
{SDE2}) and contributes to drive the drift of the other $n$ components.

Denote by $C_{0}([0,\infty)\times\tilde\Delta_{n})$ the Banach
space of
continuous functions on $[0,\infty)\times\tilde\Delta_{n}$ which vanish
at infinity, equipped with the supremum norm $\|f\|=\sup_{x\in
[0,\infty
)\times\tilde\Delta_{n}}f(x)$, and by $\mathscr{P}(B)$ the set of Borel
probability\vspace*{1pt} measures on $B$. Recall that a Markov process $\{X(t),t\ge
0\}$, taking values in a metric space $E$, is said to correspond to a
semigroup $\{T(t)\}$, acting on a closed subspace $L$ of the space of
bounded functions on $E$, if
\[
\mathbb{E}\bigl[f\bigl(X(t+s)\bigr)|\mathscr{F}_{t}^{X}\bigr]=T(s)f(X(t)),
\qquad
s,t\ge0,
\]
for every $f\in L$, where $\mathscr{F}_{t}^{X}=\sigma(X(u),u\le t)$.
%
\begin{proposition}\label{propAn-generates}
Let $\mathcal{A}_{n}$ be the operator defined in (\ref
{generatordelta-n}) and (\ref{domain-delta-n}).
The closure in $C_{0}([0,\infty)\times\tilde\Delta_{n})$ of
$\mathcal
{A}_{n}$ generates a strongly continuous, positive, conservative,
contraction semigroup $\{\mathcal{S}_{n}(t)\}$ on $C_{0}([0,\infty
)\times\tilde\Delta_{n})$. For every $\nu_{n}\in\mathscr
{P}([0,\infty
)\times\tilde\Delta_{n})$ there exists a strong Markov process
$Z^{(n)}(\cdot)=\{Z^{(n)}(t),t\ge0\}$ corresponding to $\{\mathcal
{S}_{n}(t)\}$ with initial distribution $\nu_{n}$ and sample paths in
$C_{[0,\infty)\times\tilde\Delta_{n}}([0,\infty))$ with
probability one.
\end{proposition}
\begin{pf}
We proceed by verifying the hypothesis of the Hille--Yosida theorem.
Note first that
$\mathcal{A}_{n}$ satisfies the positive maximum principle, that is,
for $f\in\mathscr{D}(\mathcal{A}_{n})$ and $(s_{*},z_{*})\in
(0,\infty
)\times\tilde\Delta_{n}$ such that $\|f\|=f(s_{*},z_{*})\ge0$ we have
$\mathcal{A}_{n}f(s_{*},z_{*})\le0$. Indeed, writing $\mathcal
{A}_{n}=\mathcal{A}_{0}+\mathcal{A}_{n,1}$
to indicate the first two and last two terms in (\ref
{generatordelta-n}), it is immediate to check that $\mathcal{A}_{0}$
and $\mathcal{A}_{n,1}$ satisfy the positive maximum principle on
$[0,\infty)$ and $\tilde\Delta_{n}$, respectively. If
$f_{0}(s_{*})\ge
0,f_{1}(z_{*})\ge0$, then $\mathcal{A}_{0}f_{0}(s)\le0$ and
$\mathcal
{A}_{n,1}f_{1}(z)\le0$, while if $f_{0}(s_{*})\le0,f_{1}(z_{*})\le0$,
then $\mathcal{A}_{0}f_{0}(s_{*})\ge0$ and $\mathcal
{A}_{n,1}f_{1}(z_{*})\ge0$. In both cases
\[
\mathcal{A}_{n}f(s_{*},z_{*})=f_{1}(z_{*})\mathcal
{A}_{0}f_{0}(s_{*})+f_{0}(s_{*})\mathcal{A}_{n,1}f_{1}(z_{*})\le0.
\]
%
Let\vspace*{1pt} now $L\subset\mathscr{D}(\mathcal{A}_{n})$ be the algebra generated
by functions $f=f_{0}\times f_{1}$, with $f_{0}\in\mathscr
{D}(\mathcal
{A}_{0})$, $\mathscr{D}(\mathcal{A}_{0})$ as in (\ref{domain-A0}) and
$f_{1}=z^{c}=z_{1}^{c_{1}}\cdots z_{n}^{c_{n}}\in C^{2}(\tilde\Delta
_{n})$, $c_{i}\in\mathbb{N}_{0}$,
so that $L$ is dense in $C_{0}([0,\infty)\times\tilde\Delta_{n})$, and
so is $\mathscr{D}(\mathcal{A}_{n})$. Denoting
$c+d_{(i)}=(c_{0},\ldots
, c_{i}+d,\ldots,c_{n})$, for $f\in L$ we have
%
\begin{eqnarray}\label{Lm->Lm-0}
&&\mathcal{A}_{n} \bigl(f_{0}(s)\times z^{ c}\bigr)\nonumber\\
&&\qquad=f_{1}(z) \mathcal{A}_{0}f(s)\nonumber\\
&&\qquad\quad{} +\frac{f_{0}(z_{0})}{2}\Biggl\{\sum_{i=1}^{n}\Biggl[ c_{i}(
c_{i}-1)\bigl(
-z^{ c}+\bigl(1-(n-2)\varepsilon_{n}\bigr)\varepsilon_{n}z^{ c-1_{(i)}}\nonumber\\
&&\hspace*{136pt}\qquad\quad{} -\varepsilon_{n}\bigl(1-(n-1)\varepsilon_{n}\bigr)z^{ c-2_{(i)}}\bigr)
\nonumber\\[-8pt]\\[-8pt]
&&\hspace*{68.8pt}\qquad\quad{} +\sum_{j\ne i}^{n}c_{i}c_{j}\bigl(-z^{c}+\varepsilon_{n}\bigl(z^{
c-1_{(j)}}+z^{ c-1_{(i)}}\nonumber\\
&&\qquad\quad\hspace*{182.5pt}{}-\varepsilon_{n}z^{ c-1_{(i)}-1_{(j)}}\bigr)
\bigr)\Biggr]
\nonumber\\
&&\hspace*{47pt}\qquad\quad{}
+\sum_{i=1}^{n} c_{i}\biggl[\frac{\beta}{z_{0}(n-1)}\bigl(z^{
c-1_{(i)}}-nz^{c}\bigr)\nonumber\\
&&\hspace*{-34pt}\hspace*{121pt}\qquad\quad{} -\alpha z^{ c-1_{(i)}}+\alpha
e^{\varepsilon_{n}\exp\{1/\varepsilon_{n}\}}e^{-\exp\{1/\varepsilon
_{n}\}z_{i}}z^{ c-\delta_{i}} \biggr]\Biggr\}\nonumber,
\end{eqnarray}
so that the image of $\mathcal{A}_{n}$ contains functions of type
$f_{0}\times z^{c}$ and $f_{0}\times e^{- b_{0}z_{i}}z^{c}$, with
$b_{0}$ fixed. For every $g(x)\in C(K)$, with $K$ compact, and
$f(x)=e^{b_{0}x}g(x)\in C(K)$, there exists a sequence $\{p^{(k)}\}$ of
polynomials on $K$ such that $\|f-p^{(k)}\|\rightarrow0$, so that
$\|e^{-b_{0}z}p^{(k)}-g\|\rightarrow0$. It follows that the image of
$\mathcal{A}_{n}$ is dense in $C_{0}([0,\infty)\times\tilde\Delta
_{n})$, and so is that of $\lambda-\mathcal{A}_{n}$ for all but at most
countably many $\lambda>0$. The first assertion now follows from
Theorem 4.2.2 of \citet{EK86} and by noting that $1\in\mathscr
{D}(\mathcal{A}_{n})$ and $\mathcal{A}_{n}1=0$, that is, $\mathcal
{A}_{n}$ is conservative.
The second assertion with $D_{[0,\infty)\times\tilde\Delta
_{n}}([0,\infty))$, the space of right-continuous functions with left
limits, in place of $C_{[0,\infty)\times\tilde\Delta_{n}}([0,\infty))$,
follows from Theorem 4.2.7 of \citet{EK86}. To prove the almost sure
continuity of sample paths, it is enough to show that for every
$z^{*}\in(0,\infty)\times\tilde\Delta_{n}$ and $\epsilon>0$ there
exists a function $f\in\mathscr{D}(\overline{\mathcal{A}}_{n})$
such that
\[
f(z^{*})=\|f\|,\qquad
\sup_{z\in B(z^{*},\epsilon)^{c}}f(z)<f(z^{*}),\qquad
\overline{\mathcal{A}}_{n}f(z^{*})=0,
\]
where $B(z^{*},\epsilon)^{c}$ is the complement of an $\epsilon
$-neighborhood of $z^{*}$ in the topology of coordinatewise convergence
[cf. \citet{EK86}, Remark~4.2.10]. This can be done by means of a
function $f\in\mathscr{D}(\mathcal{A}_{n})$ which is flat in $z^{*}$
and rapidly decreasing away from $z^{*}$, for example, of type
$f(z)=c_{1}-c_{2}\sum_{i=0}^{n}(z_{i}-z_{i}^{*})^{4}$ for appropriate
constants $c_{1},c_{2}$.
\end{pf}

For $Z^{(n)}(\cdot)$ as in Proposition \ref{propAn-generates},
consider now the mapping $\rho_{n}(Z^{(n)}(\cdot))$, where $\rho
_{n}\dvtx[0,\infty)\times\tilde\Delta_{n}\rightarrow[0,\infty)\times
\overline\nabla_{\infty}$ is defined as
%
\begin{equation}\label{rho-n}
\rho_{n}(z)=\bigl(z_{0},z_{(1)},z_{(2)},\ldots,z_{(n)},0,0,\ldots\bigr),
\end{equation}
$(z_{(1)},\ldots,z_{(n)})$ is the vector of decreasingly ordered
statistics of $(z_{1},\ldots,z_{n})\in\tilde\Delta_{n}$, and
$\overline
\nabla_{\infty}$ is the closure of the infinite-dimensional ordered
simplex, defined in (\ref{nabla-intro}).
The following proposition states that $\rho_{n}(Z^{(n)}(\cdot))$ is
still a well-defined Markov process. Define
\[
\tilde\nabla_{n}=\Biggl\{z\in\overline\nabla_{\infty}\dvtx
z_{n}\ge\varepsilon_{n}>z_{n+1}=0, \sum_{i=1}^{n}z_{i}=1
\Biggr\}
\]
and observe that $z\in\overline\nabla_{\infty}$ satisfies $z_{i}\le
1/i$ for all $i$, so that $\tilde\nabla_{n}$ is nonempty by (\ref
{varepsilon-n}).
%
\begin{proposition}\label{B-ngenerates}
Let $\tilde\mathcal{A}_{n}$ be defined by the right-hand side of
(\ref
{generatordelta-n}), with domain
\[
\mathscr{D}(\tilde\mathcal{A}_{n})
=\{f\dvtx f=f_{0}\times f_{1}, f_{0}\in\mathscr{D}(\mathcal{A}_{0}),
f_{1}\in C_{\rho_{n}}^{2}(\tilde\nabla_{n})\},
\]
where $\mathscr{D}(\mathcal{A}_{0})$ is as in (\ref{domain-A0}) and
\[
C_{\rho_{n}}^{2}(\tilde\nabla_{n})=\{f\in C(\tilde\nabla_{n})\dvtx
f\circ
\rho_{n}\in C^{2}(\tilde\Delta_{n})\}.
\]
Then the closure of $\tilde\mathcal{A}_{n}$ in $C_{0}([0,\infty
)\times
\tilde\nabla_{n})$ generates a strongly continuous, positive,
conservative, contraction semigroup $\{\mathcal{T}_{n}(t)\}$ on
$C_{0}([0,\infty)\times\tilde\nabla_{n})$. For every $\nu_{n}\in
\mathscr
{P}([0,\infty)\times\tilde\Delta_{n})$, let $Z^{(n)}$ be as in
Proposition \ref{propAn-generates}. Then $\rho_{n}(Z^{(n)}(\cdot))$ is
a strong Markov process
corresponding to $\{\mathcal{T}_{n}(t)\}$ with initial distribution
$\nu_{n}\circ\rho_{n}^{-1}$ and sample paths in $C_{[0,\infty
)\times
\tilde\nabla_{n}}([0,\infty))$ with probability one.\vspace*{-2pt}
\end{proposition}
\begin{pf}
The result follows from Proposition 2.4 in \citet{EK81}, with $\Delta
_{n}$, $\nabla_{n}$ and $\rho_{n}$ there substituted by $[0,\infty
)\times\tilde\Delta_{n}$, $[0,\infty)\times\tilde\nabla_{n}$ and
(\ref{rho-n}), respectively.\vspace*{-2pt}
\end{pf}

We now turn the attention to the limit of $\rho_{n}(Z^{(n)}(\cdot))$
when the number of types goes to infinity. To this end, consider that
$\overline\nabla_{\infty}$
is a compact and metrizable space in the topology of
coordinatewise convergence, and let
$C_{0}([0,\infty)\times\overline\nabla_{\infty})$ be the Banach space
of continuous functions on $[0,\infty)\times\overline\nabla_{\infty}$
which vanish at infinity, with the supremum norm $\|f\|=\sup_{z\in
[0,\infty)\times\overline\nabla_{\infty}}|f(z)|$.
The key issue for showing that the closure of the differential operator
$\mathcal{A}$, defined in (\ref{operatorgibbs}), generates a Feller
diffusion on $C_{0}([0,\infty)\times\overline\nabla_{\infty})$ is the
choice of the domain of $\mathcal{A}$.
Here we adapt to the present framework a technique indicated by \citet{EK81}.
Consider polynomials $\varphi_{m}\dvtx\overline\nabla_{\infty
}\rightarrow
[0,1]$ defined as\looseness=-1
%
\begin{equation}\label{power-series}
\varphi_{1}(z)=1,\qquad
\varphi_{m}(z)=\sum_{i=1}^{\infty}z_{i}^{m},\qquad
m\ge2.
\end{equation}\looseness=0
Since $z\in\overline\nabla_{\infty}$ implies $z_{i}\le i^{-1}$,
functions $\varphi_{m}$ with $m\ge2$ are uniformly convergent, and sums
in (\ref{operatorgibbs}) are assumed to be computed on
%
\begin{equation}\label{nabla}
\nabla_{\infty}=\Biggl\{z=(z_{1},z_{2},\ldots)\dvtx z_{1}\ge z_{2}\ge\cdots
\ge0, \sum_{i=1}^{\infty}z_{i}=1\Biggr\}
\end{equation}
and extended to $\overline\nabla_{\infty}$ by continuity, so that,
for example,
\[
\sum_{i=1}^{\infty}(1+z_{i})\,\frac{\partial}{\partial z_{i}}\varphi
_{2}(z)=2+2\varphi_{2}(z)
\]
instead of $2\sum_{i=1}^{\infty}z_{i}+2\varphi_{2}(z)$.
Write
%
\begin{equation}\label{gen-decomposition}
\mathcal{A}=\mathcal{A}_{0}+\mathcal{A}_{1}
\end{equation}
to indicate the first two and last two terms in (\ref{operatorgibbs}),
and denote
%
\begin{equation}\label{domain-A-1}\qquad
\mathscr{D}(\mathcal{A}_{1})
=\{\mbox{sub-algebra of }C(\overline\nabla_{\infty})
\mbox{ generated by }\varphi_{m}\mbox{ as in (\ref{power-series})}\}.\vadjust{\goodbreak}
\end{equation}
The domain $\mathscr{D}(\mathcal{A})$ of the operator (\ref
{operatorgibbs}) is then taken to be
%
\begin{eqnarray}\label{domain-B}
\mathscr{D}(\mathcal{A})
&=&\bigl\{\mbox{sub-algebra of }C_{0}\bigl([0,\infty)\times\overline
\nabla_{\infty}\bigr)\mbox{ generated by}\nonumber\\[-8pt]\\[-8pt]
&&\hspace*{38pt} f=f_{0}\times f_{1}\dvtx f_{0}\in\mathscr{D}(\mathcal{A}_{0}),
f_{1}\in\mathscr{D}(\mathcal{A}_{1})\bigr\}\nonumber
\end{eqnarray}
with $\mathscr{D}(\mathcal{A}_{0})$ as in (\ref{domain-A0}) and
$\mathscr{D}(\mathcal{A}_{1})$ as above.
%
\begin{lemma}\label{lemmaDB-dense}
The sub-algebra $\mathscr{D}(\mathcal{A}_{1})\subset C(\overline
\nabla
_{\infty})$ is dense in $C(\overline\nabla_{\infty})$.
\end{lemma}
\begin{pf}
See the proof of Theorem 2.5 in \citet{EK81}.
\end{pf}

We also need the following lemma, which shows that the operator
$\mathcal{A}_{1}$ is triangulizable.
%
\begin{lemma}\label{Lm->Lm}
Let $\mathcal{A}_{1}$ be as in (\ref{gen-decomposition}) and, for any
$m\ge2$, let $L_{m}$ be the algebra generated by polynomials as in
(\ref{power-series}), with degree not greater than $m$. Then $\mathcal
{A}_{1}\dvtx L_{m}\rightarrow L_{m}$.
\end{lemma}
\begin{pf}
The assertion follows from equation (2.4) in \citet{FS10}, with $\theta$
replaced by $\beta/s$.
\end{pf}

Then we have the following result.
%
\begin{theorem}\label{propB-generates}
Let $\mathcal{A}$ be the operator defined by (\ref{operatorgibbs})
and (\ref{domain-B}). The closure in $C_{0}([0,\infty)\times
\overline
\nabla_{\infty})$ of $\mathcal{A}$ generates a strongly continuous,
positive, conservative, contraction semigroup $\{\mathcal{T}(t)\}$ on
$C_{0}([0,\infty)\times\overline\nabla_{\infty})$.
For every $\nu\in\mathscr{P}([0,\infty)\times\overline\nabla
_{\infty
})$, there exists a strong Markov process $Z(\cdot)$ corresponding to
$\{\mathcal{T}(t)\}$ with initial distribution $\nu$ and sample paths
in $C_{[0,\infty)\times\overline\nabla_{\infty}}([0,\infty))$ with
probability one.
\end{theorem}
\begin{pf}
For every $g\in C_{0}([0,\infty)\times\overline\nabla_{\infty})$,
define $r_{n}\dvtx C_{0}([0,\infty)\times\overline\nabla_{\infty
})\rightarrow C_{0}([0,\infty)\times\tilde\nabla_{n})$ to be the
bounded linear map
\[
r_{n}g=g|_{[0,\infty)\times\tilde\nabla_{n}}
\]
given by\vspace*{1pt} the restriction of $g$ to $[0,\infty)\times\tilde\nabla_n$.
Note that $r_{n}\dvtx \mathscr{D}(\mathcal{A})\rightarrow\mathscr
{D}(\tilde\mathcal{A}_{n})$, with $\tilde\mathcal{A}_{n}$ as in Proposition \ref
{B-ngenerates}, and that
%
\begin{equation}\label{tilderng->g}
\|r_{n}g-g\|\longrightarrow0,\qquad
g\in C_{0}\bigl([0,\infty)\times\overline\nabla_{\infty}\bigr).
\end{equation}
Then, for $g\in\mathscr{D}(\mathcal{A})$ and $z\in(0,\infty)\times
\tilde\nabla_{n}$, we have
\begin{eqnarray*}
&&|\tilde\mathcal{A}_{n}r_{n}g(z)-r_{n}\mathcal{A}g(z)|\\
&&\qquad= \frac{1}{2}\Biggl|\sum_{i,j=1}^{n}
\bigl(a_{ij}^{(n)}(z)-z_{i}(\delta_{ij}-z_{j})\bigr)\,\frac{\partial^{2}
g(z)}{\partial z_{i}\,\partial z_{j}}\\
&&\hspace*{8.6pt}\qquad\quad{} +\sum_{i=1}^{n}\biggl(\frac{\beta}{z_{0}(n-1)}(1-z_{i})-\alpha
\bigl(\exp\{-(z_{i}-\varepsilon_{n})e^{1/\varepsilon_{n}}\}\bigr)
\biggr)\,\frac{\partial g(z)}{\partial z_{i}}\Biggr|
\end{eqnarray*}
with $a_{ij}^{(n)}(z)$ as in (\ref{a-ij}).
In particular,
\begin{eqnarray*}
&&
\bigl|a_{ij}^{(n)}(z)-z_{i}(\delta_{ij}-z_{j})\bigr|\\
&&\qquad=
\cases{\displaystyle
\varepsilon_{n}[(z_{i}-\varepsilon_{n})(n-1)+1-z_{i}],&\quad
$i=j$,\vspace*{2pt}\cr
\displaystyle
\varepsilon_{n}[z_{i}+z_{j}-\varepsilon_{n}],&\quad $i\ne j$,}
\end{eqnarray*}
which is bounded above by $n\varepsilon_{n}$, from which
%
\begin{eqnarray}\label{majoration}
&&
|\tilde\mathcal{A}_{n}r_{n}g(z)-r_{n}\mathcal{A}g(z)|\nonumber\\
&&\qquad\le n\varepsilon_{n}\sum_{i,j=1}^{n}\biggl|\frac{\partial^{2}
g(z)}{\partial z_{i}\,\partial z_{j}}\biggr| +\frac{\beta}{z_{0}(n-1)}\sum_{i=1}^{n}\biggl|\frac{\partial
g(z)}{\partial z_{i}}\biggr|\\
&&\qquad\quad{}
+\sum_{i=1}^{n}\exp\{-(z_{i}-\varepsilon_{n})e^{1/\varepsilon
_{n}}\}\biggl|\frac{\partial g(z)}{\partial z_{i}}\biggr|\nonumber.
\end{eqnarray}
For $g\in\mathscr{D}(\mathcal{A})$ of type $g=f_{0}\times f_{1}$, with
$f_{0}\in\mathscr{D}(\mathcal{A}_{0})$ and $f_{1}=\varphi
_{m_{1}}\cdots
\varphi_{m_{k}}$, we have
%
\begin{eqnarray}\label{firstderivative}
\sum_{i=1}^{n}\biggl|\frac{\partial g(z)}{\partial z_{i}}\biggr|
&=&|f_{0}(z_{0})|\sum_{i=1}^{n}\sum_{j=1}^{k}m_{j}z_{i}^{m_{j}-1}\prod
_{h\ne j}\varphi_{m_{h}}\nonumber\\[-8pt]\\[-8pt]
&\le&|f_{0}(z_{0})|\sum_{j=1}^{k}m_{j}\sum_{i=1}^{n}z_{i}^{m_{j}-1}\nonumber
\end{eqnarray}
so that
\[
\sum_{i=1}^{n}\exp\{-(z_{i}-\varepsilon_{n})e^{1/\varepsilon
_{n}}\}\biggl|\frac{\partial g(z)}{\partial z_{i}}\biggr|
\le
n\varepsilon_{n}|f_{0}(z_{0})|\sum_{j=1}^{k}m_{j}\rightarrow0
\]
%
uniformly as $n\rightarrow\infty$ by (\ref{varepsilon-n}).
Furthermore,
\begin{eqnarray*}
&&
\sum_{i,j=1}^{n} \biggl|\frac{\partial^{2} g(z)}{\partial
z_{i}\,\partial z_{j}}\biggr|\\
&&\qquad\le|f_{0}(z_{0})|\sum_{i,j=1}^{\infty}\biggl[\partial_{ij}\varphi_{m_{h}}
\prod_{\ell\ne h} \varphi_{m_{\ell}}+\sum_{q\ne h}\partial
_{i}\varphi
_{m_{h}}\partial_{j}\varphi_{m_{q}}\prod_{\ell\ne h,q} \varphi
_{m_{\ell
}}\biggr]\\
&&\qquad= |f_{0}(z_{0})|\biggl[m_{h}(m_{h}-1)\varphi_{m_{h}-2}
\prod_{\ell\ne h} \varphi_{m_{\ell}}\\
&&\hspace*{38pt}\qquad\quad{}+\sum_{q\ne h}m_{h}m_{q}\varphi_{m_{h}+m_{q}-2}\prod_{\ell\ne h,q}
\varphi_{m_{\ell}}\\
&&\hspace*{38pt}\qquad\quad{} +\sum_{q\ne h}m_{h}m_{q}\varphi_{m_{h}-1}\varphi_{m_{q}-1}\prod
_{\ell
\ne h,q} \varphi_{m_{\ell}}
\biggr]\\
&&\qquad\le
|f_{0}(z_{0})|\biggl[m_{h}(m_{h}-1)+2\sum_{q\ne h}m_{h}m_{q}
\biggr],
\end{eqnarray*}
whose right-hand side is bounded. Since also the right-hand side of
(\ref{firstderivative}) is bounded above by $|f_{0}(z_{0})|\sum
_{j=1}^{k}m_{j}$,
it follows by (\ref{varepsilon-n}) that the right-hand side of (\ref
{majoration}) goes to zero uniformly and, by means of (\ref
{tilderng->g}), that
%
\begin{equation}\label{eqBn-convergence}
\|\tilde\mathcal{A}_{n}r_{n}g-\mathcal{A}g\|\longrightarrow0,\qquad
g\in\mathscr{D}(\mathcal{A}).
\end{equation}
Proposition \ref{B-ngenerates} implies that $\tilde\mathcal{A}_{n}$ is
a dissipative operator for every $n\ge1$, so that by (\ref
{eqBn-convergence}) $\mathcal{A}$ is dissipative. Moreover, Lemma
\ref{lemmaDB-dense} and Lemma \ref{Lm->Lm}, respectively, imply that
$\mathscr{D}(\mathcal{A})$ and the range of $\lambda-\mathcal{A}$, for
all but at most countably many $\lambda>0$, are dense in
$C_{0}([0,\infty)\times\overline\nabla_{\infty})$. The fact that the
closure of $\mathcal{A}$ generates a strongly continuous contraction
semigroup $\{\mathcal{T}(t)\}$ on
$C_{0}([0,\infty)\times\overline\nabla_{\infty})$ now follows from the
Hille--Yosida theorem [see Theorem~1.2.12 in \citet{EK86}]. It is
also immediate to check that $\mathcal{A}1=0$, so that $(1,0)\in
\overline{(g,\mathcal{A}g)}$ and $\{\mathcal{T}(t)\}$ is conservative.
Finally, (\ref{eqBn-convergence}), together with Lemma
\ref{lemmaDB-dense} and Theorem 1.6.1 of \citet{EK86}, implies the
semigroup convergence
%
\begin{equation}\label{eqSG-convergence}
\|\mathcal{T}_{n}(t)r_{n}g-\mathcal{T}(t)g\|\longrightarrow0,\qquad
g\in C_{0}\bigl([0,\infty)\times\overline\nabla_{\infty}\bigr),
\end{equation}
uniformly on bounded intervals. From Proposition \ref{B-ngenerates},
$\{\mathcal{T}_{n}(t)\}$ is a positive operator for every $n\ge1$, so
that $\{\mathcal{T}(t)\}$ is in turn positive.

The second\vspace*{1pt} assertion of the theorem, with $D_{[0,\infty)\times
\overline
\nabla_{\infty}}([0,\infty))$ in place of $C_{[0,\infty)\times
\overline
\nabla_{\infty}}([0,\infty))$, follows\vspace*{1pt} from Theorem 4.2.7 in \citet
{EK86}, while the continuity of sample paths follows from a similar
argument to that used in the proof of Proposition \ref{propAn-generates}.
\end{pf}

The following corollary formalizes the convergence in distribution of the sequence
of processes of Proposition \ref{B-ngenerates} to the
infinite-dimensional diffusion of Theorem \ref{propB-generates}.
%
\begin{corollary}
Let $Z^{(n)}(\cdot)$ be as in Proposition \ref{propAn-generates} with
initial distribution $\nu_{n}\in\mathscr{P}([0,\infty)\times\tilde
\Delta
_n)$, and let $Z(\cdot)$ be as in Theorem \ref{propB-generates} with
initial distribution $\nu\in\mathscr{P}([0,\infty)\times\overline
\nabla
_{\infty})$. If $\nu_{n}\circ\rho_{n}^{-1}\Rightarrow\nu$ on
$\overline
\nabla_{\infty}$, then
\[
\rho_{n}\bigl(Z^{(n)}(\cdot)\bigr)\Rightarrow Z(\cdot)
\qquad\mbox{in }C_{[0,\infty)\times\overline\nabla_{\infty}}([0,\infty))
\]
as $n\rightarrow\infty$.
\end{corollary}
\begin{pf}
The result\vspace*{1pt} with $D_{[0,\infty)\times\overline\nabla_{\infty
}}([0,\infty
))$ in place of $C_{[0,\infty)\times\overline\nabla_{\infty
}}([0,\infty
))$ follows from Proposition \ref{B-ngenerates}, together with (\ref
{eqSG-convergence}) and Theorem 4.2.5 in \citet{EK86}. The fact that
the weak convergence holds in $C_{[0,\infty)\times\overline\nabla
_{\infty}}([0,\infty))$ follows from relativization of the Skorohod topology.
\end{pf}


\section{A population model for normalized inverse-Gaussian
diffusions}\label{sectionpopulation}

By formalizing the population process briefly mentioned in Section \ref
{sections-diffusion} for constructing the $\alpha$-diversity
diffusion, in this section we provide a discrete approximation, based
on a countable number of particles, for the diffusion with operator
(\ref{operatorgibbs}). More specifically, this is obtained as the
limit in distribution of the process of frequencies of types associated
with a set of particles sampled from a normalized inverse-Gaussian
random probability measure, jointly with the normalized version of the
diversity process.

In view of (\ref{Pitman-urn}), the conditional distribution of the
$i$th component of an exchangeable sequence $(X_{1},\ldots,X_{n})$
drawn from a random probability measure of Gibbs type can be written
%
\begin{eqnarray}\label{eqi-predictive}
\mathbb{P}\{X_{i}&\in&\cdot| X_{1},\ldots,X_{i-1},X_{i+1},\ldots
,X_{n}\}\nonumber\\
&=& g_{0}(n-1,K_{n-1,i})\nu_0(\cdot)\\
&&{}+g_{1}(n-1,K_{n-1,i})\sum
_{j=1}^{K_{n-1,i}}(n_j-\alpha) \delta_{X_j^*}(\cdot),\nonumber
\end{eqnarray}
where $\nu_{0}$ is a nonatomic probability measure and
$(X_{1}^{*},\ldots,X_{K_{n-1,i}}^{*})$ are the $K_{n-1,i}$ distinct
values in $(X_{1},\ldots,X_{i-1},X_{i+1},\ldots,X_{n})$. For fixed $n$,
define a Markov chain $\{X^{(n)}(m),m\ge0\}$ on $\mathbb{X}^{n}$ by
means of the transition semigroup
\[
T_{n}f(x)=\int f(y)p_{n}(x,\dd y),\qquad
f\in C_{0}(\mathbb{X}^{n}),
\]
where $x,y\in\mathbb{X}^{n}$, $C_{0}(\mathbb{X}^{n})$ is the space of
Borel-measurable continuous functions on $\mathbb{X}^{n}$ vanishing at
infinity,
%
\begin{equation}\label{particle-transition}
p_{n}(x,\dd y)=\frac{1}{n}\sum_{i=1}^{n}\tilde p_{1}\bigl(\dd
y_{i}|x_{(-i)}\bigr)\prod_{k\ne i}\delta_{x_{k}}(\dd y_{k}),
\end{equation}
$x_{(-i)}=(x_{1},\ldots,x_{i-1},x_{i+1},\ldots,x_{n})$ and $\tilde
p_{1}(\dd y|x_{(-i)})$ is (\ref{eqi-predictive}). The interpretation is
as follows. At each transition one component is selected at random with
uniform probability, and is updated with a value sampled from (\ref
{eqi-predictive}), conditional on all other components, which are left
unchanged. Hence, the incoming particle is either a new type (a mutant
offspring) or a copy of an old type (a copied offspring). Embed now the
chain in a pure jump Markov process on $\mathbb{X}^{n}$ with
exponentially distributed waiting times with intensity one, and denote
the resulting process by $X^{(n)}(\cdot)=\{X^{(n)}(t),t\ge0\}$.\vadjust{\goodbreak}
The infinitesimal generator of $X^{(n)}(\cdot)$ is given by
%
\begin{eqnarray}\label{particle-gen}
B_{n}f(x)
&=&
\frac{1}{n}\sum_{i=1}^ng_{0}^{(n-1,i)}[P_{i}f(x)-f(x)]\nonumber\\[-8pt]\\[-8pt]
&&{} +\frac{1}{n}\sum_{i=1}^n\sum
_{j=1}^{K_{n-1,i}}g_{1}^{(n-1,i)}(n_{j}-\alpha)[\Phi
_{j^{*}i}f(x)-f(x)]\nonumber
\end{eqnarray}
with domain
%
\begin{equation}\label{particle-domain}
\mathscr{D}(B_{n})=\{f\dvtx f\in C_{0}(\mathbb{X}^{n})\}.
\end{equation}
Here $\Phi_{j^{*}i}\dvtx C_{0}(\mathbb{X}^{n})\rightarrow C_{0}(\mathbb
{X}^{n-1})$ is defined as
%
\begin{equation}\label{Phi}
\Phi_{j^{*}i}f(x_{1},\ldots,x_{n})=f(x_{1},\ldots
,x_{i-1},x_{j}^{*},x_{i+1},\ldots,x_{n})
\end{equation}
for $x_{j}^{*}\in(x_{1}^{*},\ldots,x_{K_{n-1,i}}^{*})$, $P$ is the
transition semigroup
%
\begin{equation}\label{mutation}
Pg(z)=\int g(y)p_{1}(z,\dd y),\qquad
g\in C_{0}(\mathbb{X}),
\end{equation}
where $p_{1}(z,\dd y)$ is given by
%
\begin{equation}\label{nu0}
p_{1}(z,\dd y)=\nu_{0}(\dd y),
\end{equation}
$P_{i}f$ denotes $P$ acting on the $i$th coordinate of $f$, and we have
set for brevity
%
\begin{equation}\label{pesi-short}
g_{j}^{(n-1,i)}=g_{j}(n-1,K_{n-1,i}),\qquad
j=0,1.
\end{equation}
%
Defining (\ref{mutation}) and (\ref{nu0}) separately is somewhat
redundant, but will allow us to provide a general expression for the
global mutation rate in this particle representation before making the
assumptions of nonatomicity and parent independence as in (\ref{nu0}).
See (\ref{mutationrate}) below.

Define now the map $w\dvtx\mathbb{X}^{n}\rightarrow\overline\nabla
_{\infty
}$ by
%
\begin{equation}\label{gamma}
w(x)=w\bigl(x^{(n)}\bigr)=(z_{1},\ldots,z_{K_{n}},0,0,\ldots),
\end{equation}
where $z_{j}$ and $K_{n}$, respectively, denote the relative frequency
of the $j$th most abundant type and the number of types in $X^{(n)}$.
Let also $\mathcal{A}$ be as in (\ref{operatorgibbs}). The next
theorem states that
%
\begin{equation}\label{[K,Z]}
\bigl[K_{n}(\cdot)/n^{\alpha},w\bigl(X^{(n)}(\cdot)\bigr)\bigr]
=\bigl\{\bigl[K_{n}(t)/n^{\alpha},w\bigl(X^{(n)}(t)\bigr)\bigr],t\ge0\bigr\},
\end{equation}
if appropriately\vspace*{1pt} rescaled in time, converges in
distribution to the process with generator $\mathcal{\overline A}$. The
proof is deferred to the \hyperref[app]{Appendix} and contains, as a byproduct, a~more
formal derivation of Theorem \ref{S-convergence}.
%
\begin{theorem}\label{theoremparticle-convergence}
Let $X^{(n)}(\cdot)$ be the $\mathbb{X}^{n}$-valued process with
generator (\ref{particle-gen}) and (\ref{particle-domain}), $w\dvtx\mathbb
{X}^{n}\rightarrow\overline\nabla_{\infty}$ as in (\ref{gamma}) and
$Z(\cdot)$ as in Theorem \ref{propB-generates}. If
%
\begin{equation}\label{conv-distr-iniz}
\bigl[K_{n}(0)/n^{\alpha},w\bigl(X^{(n)}(0)\bigr)\bigr]\Rightarrow Z(0),
\end{equation}
then
%
\begin{equation}\label{eqparticle-convergence}
\bigl[K_{n}(n^{3/2}t)/n^{\alpha},w\bigl(X^{(n)}(n^{2}t/2)\bigr)\bigr]\Rightarrow Z(t)
\end{equation}
in $C_{[0,\infty)\times\overline\nabla_{\infty}}([0,\infty))$.
\end{theorem}

We conclude the section by showing the reversibility of the particle process.
Denote the joint distribution of an $n$-sized sequence from the
generalized P\'olya urn scheme (\ref{Pitman-urn}) by
\begin{eqnarray*}
&&\mathscr{M}_{n}(\dd x_{1},\ldots,\dd x_{n})\\
&&\qquad= \nu_{0}(\dd x_{1})\prod_{i=1}^{n-1}
\Biggl[g_{0}(i,K_{i})\nu_0(\dd x_{i+1})
+g_{1}(i,K_{i})\sum_{j=1}^{K_{i}}(n_j-\alpha) \delta_{x_j^*}(\dd
x_{i+1})\Biggr].
\end{eqnarray*}

\begin{proposition}
Let $X^{(n)}(\cdot)$ be the $\mathbb{X}^{n}$-valued process with
generator given by (\ref{particle-gen}) and (\ref{particle-domain}). Then
$X^{(n)}(\cdot)$ is reversible with respect to $\mathscr{M}_{n}$.
\end{proposition}
\begin{pf}
Let $q_{n}(x,\dd y)$ denote the infinitesimal transition kernel on
$\mathbb{X}^{n}\times\mathscr{B}(\mathbb{X}^{n})$ of $X^{(n)}$.
%
Denoting by $\lambda_{n}$ the rate at which the discontinuities
of $X^{(n)}$ occur, and
recalling (\ref
{particle-transition}), we have
\begin{eqnarray*}
&&\mathscr{M}_{n}(\dd x)q_{n}(x,\dd y)\\
&&\qquad= \mathscr{M}_{n}(\dd x)\lambda_{n}
\frac{1}{n}\sum_{i=1}^{n}p_{1}\bigl(\dd y_{i}|x_{(-i)}\bigr)\prod_{k\ne
i}\delta
_{x_{k}}(\dd y_{k})\\
&&\qquad= \frac{\lambda_{n}}{n}\sum_{i=1}^{n}\mathscr{M}_{n-1}\bigl(\dd x_{(-i)}\bigr)p_{1}
\bigl(\dd x_{i}|x_{(-i)}\bigr)p_{1}\bigl(\dd y_{i}|x_{(-i)}\bigr)\prod_{k\ne i}\delta
_{x_{k}}(y_{k})\nonumber\\
&&\qquad= \frac{\lambda_{n}}{n}\sum_{i=1}^{n}\mathscr{M}_{n-1}\bigl(\dd y_{(-i)}\bigr)p_{1}
\bigl(\dd x_{i}|y_{(-i)}\bigr)p_{1}\bigl(\dd y_{i}|y_{(-i)}\bigr)\prod_{k\ne i}\delta
_{y_{k}}(x_{k})\nonumber\\
&&\qquad= \mathscr{M}_{n}(\dd y)\frac{1}{n}\sum_{i=1}^{n}\lambda_{n}p_{n}
(\dd x_{i}|y_{-i})\prod_{k\ne i}\delta_{y_{k}}(x_{k})
=\mathscr{M}_{n}(\dd y)q_{n}(y,\dd x),
\end{eqnarray*}
giving the result.
\end{pf}


\section{Conditioning on the alpha-diversity}\label{secconditioning}

We conclude by discussing an interesting connection with the
two-parameter model (\ref{operatortheta-sigma}). In the \hyperref[intro]{Introduction}
it was observed that conditioning on the $\alpha$-diversity diffusion\vadjust{\goodbreak}
$S_{t}$ to be constant in the operator (\ref{operatorgibbs}) only
yields the special case $\beta=0$ and $\alpha=1/2$, consistently with
the associated random probability measures. It turns out that
performing the same conditioning operation in the particle construction
of the previous section, before taking the limit for $n\rightarrow
\infty
$, yields a particular instance of the two-parameter model.
The following proposition states that under this pre-limit conditioning
with $S_{t}\equiv s$, the normalized inverse-Gaussian model with
operator (\ref{operatorgibbs}) reduces to the two-parameter model
with $(\theta,\alpha)=((\alpha s)^{2},\alpha)$ and $\alpha=1/2$.
%
\begin{proposition}
Let $X^{(n)}(\cdot)$ be as in Theorem \ref{theoremparticle-convergence}, $w$
be as in (\ref{gamma}), $\tilde Z^{(n)}(\cdot)$~be defined by the
left-hand side of (\ref{eqparticle-convergence}),
and denote by $V^{\theta,\alpha}(\cdot)$ the process with operator
$\mathcal{L}^{\theta,\alpha}$ as in (\ref{operatortheta-sigma}). Then
\[
\bigl[\tilde Z^{(n)}(\cdot)| S_{t}\equiv s\bigr]\Rightarrow V^{s^{2}/4,1/2}(\cdot)
\]
in $C_{[0,\infty)\times\overline\nabla_{\infty}}([0,\infty))$ as
$n\rightarrow\infty$.
\end{proposition}
\begin{pf}
Let $\alpha=1/2$ throughout the proof.
In the pre-limit version of the process of frequencies derived from the
particle process, that is, (\ref{gamma}), conditioning on $S_{t}\equiv
s$ means conditioning on $K_{n}(\cdot)$ being constant over time, hence
with zero conditional first and second moment.
Denote
\[
z-\varepsilon_{i}=(z_{1},\ldots,z_{i-1},z_{i}-1,z_{i+1},\ldots)
\]
and assume $z$ has $k$ nonnull components obtained from $n$ particles.
Then, as in Section \ref{sectionpopulation}, when a particle is
removed we have the change of frequency
\[
z\mapsto z-\frac{\varepsilon_{i}}{n} \qquad\mbox{w.p. }z_{i},
\]
where $z-\varepsilon_{i}/n$ has
\begin{eqnarray*}
&&\mbox{(1):}\quad k\mbox{ nonnull components} \qquad\mbox{w.p. }
1-\frac{m_{1,n}}{n}, \\
&&\mbox{(2):}\quad k-1\mbox{ nonnull components} \qquad\mbox{w.p. }
\frac{m_{1,n}}{n}.
\end{eqnarray*}
Conditional on case (1), the number of nonnull components remains $k$
if the incoming particle is a copy of an existing type, that is, we
observe either of
\begin{eqnarray*}
z-\frac{\varepsilon_{i}}{n}&\mapsto& z-\frac{\varepsilon
_{i}}{n}+\frac{\varepsilon_{i}}{n}
\qquad\mbox{w.p. } g_{1}^{(n,k)}(n_{i}-1-\alpha
)/\bigl(1-g_{0}^{(n,k)}\bigr),\\
z-\frac{\varepsilon_{i}}{n}&\mapsto& z-\frac{\varepsilon
_{i}}{n}+\frac{\varepsilon_{j}}{n}
\qquad\mbox{w.p. } g_{1}^{(n,k)}(n_{j}-\alpha)/\bigl(1-g_{0}^{(n,k)}\bigr),
\end{eqnarray*}
%
where $g_{0}^{(n,k)}$ and $g_{1}^{(n,k)}$ are as in (\ref{pesi-short}),
while conditional on case (2) we observe
\[
z-\frac{\varepsilon_{i}}{n}\mapsto z-\frac{\varepsilon
_{i}}{n}+\frac{\varepsilon_{k+1}}{n}
\qquad\mbox{w.p. }1.
\]
%
For $\lambda_{n}=n^2/2$,
the generator of the process
$\tilde Z^{(n)}(\cdot)$ in this case can be written
\begin{eqnarray*}
\hspace*{-4pt}&&\mathcal{B}_{n,1}f_{1}(z)\\
\hspace*{-4pt}&&\qquad= \lim_{\delta t\downarrow0}\frac{1}{\delta t}
\Biggl\{\lambda_{n}\delta t
\sum_{i=1}^{k}z_{i}\biggl[
f_{1}(z) \biggl(1-\frac{m_{1,n}}{n}\biggr)\frac
{g_{1}^{(n,k)}(n_{i}-1-\alpha)}{1-g_{0}^{(n,k)}}\\
\hspace*{-4pt}&&\hspace*{86pt}\qquad\quad{} +\sum_{j\ne i}f_{1}\biggl(z-\frac{\varepsilon_{i}}{n}+\frac
{\varepsilon_{j}}{n}\biggr) \biggl(1-\frac{m_{1,n}}{n}\biggr)\frac
{g_{1}^{(n,k)}(n_{j}-\alpha)}{1-g_{0}^{(n,k)}}\\
\hspace*{-4pt}&&\hspace*{190pt}\qquad\quad{} + f_{1}\biggl(z-\frac{\varepsilon_{i}}{n}+\frac{\varepsilon
_{k+1}}{n}\biggr) \frac{m_{1,n}}{n}\biggr]\\
\hspace*{-4pt}&&\hspace*{146.2pt}\qquad\quad{} +(1- \lambda_{n}\delta t)f_{1}(z)+O((\delta t)^{2})-f_{1}(z)\Biggr\}
\end{eqnarray*}
for $f_{1}\in C^{2}(\nabla_{n})$,
%
\begin{equation}\label{nabla-n}
\nabla_{n}=\Biggl\{z\in\nabla_{\infty}\dvtx z_{n+1}=0, \sum
_{i=1}^{n}z_{i}= 1\Biggr\},
\end{equation}
and $\nabla_{\infty}$ as in (\ref{nabla}).
Exploiting the relation
\[
\biggl(1-\frac{m_{1,n}}{n}\biggr)\frac{g_{1}\sum_{j=1}^{k}(n_{j}-\alpha
)}{1-g_{0}}+\frac{m_{1,n}}{n}=1,
\]
we can write
%
\begin{eqnarray}\label{a-n,1}
&&\mathcal{B}_{n,1}f_{1}(z)\nonumber\\
&&\qquad= \lambda_{n}
\sum_{i=1}^{k}z_{i}\biggl\{
[f_{1}(z)-f_{1}(z)] \biggl(1-\frac{m_{1,n}}{n}\biggr)\frac
{g_{1}^{(n,k)}(n_{i}-1-\alpha)}{1-g_{0}^{(n,k)}}\nonumber\\
&&\hspace*{40.8pt}\qquad\quad{} +\sum_{j\ne i}
\biggl[f_{1}\biggl(z-\frac{\varepsilon_{i}}{n}+\frac{\varepsilon
_{j}}{n}\biggr)-f_{1}(z)\biggr] \biggl(1-\frac{m_{1,n}}{n}\biggr)\frac
{g_{1}^{(n,k)}(n_{j}-\alpha)}{1-g_{0}^{(n,k)}}\nonumber\\[-8pt]\\[-8pt]
&&\hspace*{142.5pt}\qquad\quad{} + \biggl[f_{1}\biggl(z-\frac{\varepsilon_{i}}{n}+\frac{\varepsilon
_{k+1}}{n}\biggr)-f_{1}(z)\biggr]\frac{m_{1,n}}{n}\biggr\}\nonumber\\
&&\qquad= \lambda_{n}\sum_{i,j}z_{i}
\biggl[f_{1}\biggl(z-\frac{\varepsilon_{i}}{n}+\frac{\varepsilon
_{j}}{n}\biggr)-f_{1}(z)\biggr] \biggl(1-\frac{m_{1,n}}{n}\biggr)\frac
{g_{1}^{(n,k)}(n_{j}-\alpha)}{1-g_{0}^{(n,k)}}\nonumber\\
&&\qquad\quad{} + \lambda_{n}\sum_{i=1}^{k}z_{i}\biggl[f_{1}\biggl(z-\frac{\varepsilon
_{i}}{n}+\frac{\varepsilon_{k+1}}{n}\biggr)-f_{1}(z)\biggr] \frac
{m_{1,n}}{n}.\nonumber
\end{eqnarray}
By making use of Taylor's theorem, it can be easily verified that the
following three relations holds:
\begin{eqnarray*}
\sum_{i,j}z_{i}z_{j}
\biggl[f_{1}\biggl(z-\frac{\varepsilon_{i}}{n}+\frac{\varepsilon
_{j}}{n}\biggr)-f_{1}(z)\biggr]
&=& \frac{1}{n^{2}}\sum_{i,j}z_{i}(\delta_{ij}-z_{j})\,\frac{\partial
^{2}f_{1}}{\partial z_{i}\,\partial z_{j}}+o(n^{-2}),\\
\sum_{i,j}z_{i}
\biggl[f_{1}\biggl(z-\frac{\varepsilon_{i}}{n}+\frac{\varepsilon
_{j}}{n}\biggr)-f_{1}(z)\biggr]
&=& \frac{1}{n}\sum_{i}\frac{\partial f_{1}}{\partial z_{i}}-\frac
{k}{n}\sum_{i}z_{i}\,\frac{\partial f_{1}}{\partial z_{i}}+O(n^{-2}),\\
\sum_{i}z_{i}
\biggl[f_{1}\biggl(z-\frac{\varepsilon_{i}}{n}+\frac{\varepsilon
_{k+1}}{n}\biggr)-f_{1}(z)\biggr]
&=& -\frac{1}{n}\sum_{i}z_{i}\,\frac{\partial f_{1}}{\partial z_{i}}+O(n^{-2}).
\end{eqnarray*}
By means of the last three expressions, we can write (\ref{a-n,1}) as
\[
\mathcal{B}_{n,1}f_{1}(z)= \frac{\lambda_{n}}{n^{2}}\sum
_{i,j}z_{i}(\delta_{ij}-z_{j})\,\frac{\partial^{2}f_{1}(z)}{\partial
z_{i}\,\partial z_{j}}
+\mathcal{B}_{n,1}^{\mathrm{dr}}f_{1}(z)+o(\lambda_{n}n^{-2}),
\]
where $\mathcal{B}_{n,1}^{\mathrm{dr}}f_{1}(z)$ is the drift term,
given by
\begin{eqnarray*}
\mathcal{B}_{n,1}^{\mathrm{dr}}f_{1}(z)&=&
\frac{\lambda_{n}}{n}\biggl[\biggl(\frac{\alpha k g_{1}}{1-g_{0}}
\biggl(1-\frac{m_{1,n}}{n}\biggr)-\frac{m_{1,n}}{n}\biggr)\sum_{i}z_{i}\,\frac
{\partial f_{1}}{\partial z_{i}}\\
&&\hspace*{65pt}{} -\frac{\alpha g_{1}}{1-g_{0}}\biggl(1-\frac{m_{1,n}}{n}\biggr)\sum
_{i}\,\frac{\partial f_{1}}{\partial z_{i}}\biggr].
\end{eqnarray*}
Using (\ref{asympdistribution}) and Proposition \ref
{second-order-approximation},
%
%
it can be seen that
\[
\frac{\alpha k g_{1}}{1-g_{0}}\biggl(1-\frac{m_{1,n}}{n}\biggr)-\frac{m_{1,n}}{n}
\approx-\frac{\alpha k m_{1,n}}{n^{2}}
\approx-\frac{\alpha^{2}s^{2}}{n}
\]
and
\[
\frac{\alpha g_{1}}{1-g_{0}}\biggl(1-\frac{m_{1,n}}{n}\biggr)
\approx\frac{\alpha}{n},
\]
yielding
\begin{eqnarray*}
\mathcal{B}_{n,1}f_{1}(z)&=&
\frac{\lambda_{n}}{n^{2}}\sum_{i,j}z_{i}(\delta_{ij}-z_{j})\,\frac
{\partial^{2}f_{1}(z)}{\partial z_{i}\,\partial z_{j}}\\
&&{} -\frac{\lambda_{n}}{n^{2}}\sum_{i}\bigl((\alpha s)^{2}z_{i}+\alpha
\bigr)\,\frac
{\partial f_{1}(z)}{\partial z_{i}}+o(\lambda_{n}n^{-2}).
\end{eqnarray*}
For $f\in C(\overline\nabla_{\infty})$, define $\tilde
r_{n}\dvtx C(\overline
\nabla_{\infty})\rightarrow C(\nabla_{n})$ to be
%
\begin{equation}\label{rn-tilde}
\tilde r_{n}f=f|_{\nabla_{n}},
\end{equation}
namely, the restriction of $f$ to $\tilde\nabla_{n}$, and note that
\[
\|\tilde r_{n}f-f\|\longrightarrow0,\qquad
f\in C(\overline\nabla_{\infty}).
\]
Recalling that $\lambda_{n}=n^{2}/2$ implies that for $f$ as in (\ref
{domain-A-1}) and $\mathcal{L}^{\theta,\alpha}$ as in (\ref
{operatortheta-sigma}),
\[
\bigl\|\mathcal{L}^{(\alpha s)^{2},\alpha}f-\mathcal{B}_{n,1}\tilde
r_{n}f\bigr\|\rightarrow0.
\]
The strong convergence of the corresponding semigroups on $C(\overline
\nabla_{\infty})$, similar to (\ref{eqSG-convergence}), and the
statement of the proposition now follow from an application of Theorems
1.6.1 and 4.2.11 in \citet{EK86}, together with the relativization of
the Skorohod topology to $C_{[0,\infty)\times\overline\nabla
_{\infty
}}([0,\infty))$.
\end{pf}


\begin{appendix}\label{app}
\section*{Appendix}

\subsection*{Proof of Proposition \protect\ref{second-order-approximation}}%
Consider first that $V_{n,k}$ appearing in (\ref{g0-g1}), in the case
of generalized gamma processes, can be written [cf. \citet{LMP07b}]
\[
V_{n,k}=\frac{a^{k}}{\Gamma(n)}\int_{0}^{\infty}x^{n}\exp\biggl\{-\frac
{a}{\alpha}[(\tau+x)^{\alpha}-\tau^{\alpha}]\biggr\}(\tau+x)^{\alpha
k-n}\,\dd x.
\]
Together with (\ref{V-n,k-recursion}), this leads to writing
%
\setcounter{equation}{0}
\begin{equation}\label{decompos-g0}
g_{0}(n,k)
=\frac{V_{n,k}-(n-\alpha k)V_{n+1,k}}{V_{n,k}}
=1-(1-\alpha k/n)w(n,k),
\end{equation}
where
\[
w(n,k)=\frac{\int_{0}^{\infty}x^{n}\exp\{-({a}/{\alpha})[(\tau
+x)^{\alpha}-\tau^{\alpha}]\}(\tau+x)^{\alpha k-n-1}\,\dd x}
{\int_{0}^{\infty}x^{n-1}\exp\{-({a}/{\alpha})[(\tau+x)^{\alpha
}-\tau
^{\alpha}]\}(\tau+x)^{\alpha k-n}\,\dd x}.
\]
Denote by $f(x)$ the integrand of the denominator of $w(n,k)$, so
\[
w(n,k)=\int_{0}^{\infty}\frac{x}{\tau+x}f(x)\,\dd x\Big/\int_{0}^{\infty
}f(x)\,\dd x.
\]
Since $f(x)$ is unimodal, by means of the Laplace method one can
approximate $f(x)$ with the kernel of a normal density with mean given by
%
\begin{equation}\label{argmax}
x^{*}=\argmax_{x>0} x^{n-1}\exp\biggl\{-\frac{a}{\alpha}[(\tau
+x)^{\alpha
}-\tau^{\alpha}]\biggr\}(\tau+x)^{\alpha k-n}
\end{equation}
and variance given by $-[f''(x)]^{-1}\vert_{x=x^{*}}$.
It follows that
\[
w(n,k)\approx\frac{f_{N}(x^{*}_{N})C(x^{*}_{N},-[f''_{N}(x)]^{-1}
\vert_{x=x^{*}})}
{f(x^{*}_{D})C(x^{*}_{D},-[f''(x)]^{-1}\vert_{x=x^{*}})},
\]
where $f_{N}$ denotes the integrand of the numerator, $x^{*}_{N}$ and
$x^{*}_{D}$ the modes of the integrands of numerator and denominator,
respectively, and $C(x,y)$ is the normalizing constant of a normal
kernel with mean $x$ and variance $y$, yielding
%
\begin{equation}\label{approx-phi}
w(n,k)\approx
\frac{f_{N}(x^{*}_{N})}{f(x^{*}_{D})}
\biggl(\frac{f''(x^{*}_{D})}{f''_{N}(x^{*}_{N})}\biggr)^{1/2}.
\end{equation}
From (\ref{argmax}), the mode $x^{*}_{D}$ is the only positive real
root of the equation
%
\begin{equation}\label{polynomial}
(n-1)x^{-1}+(\alpha k-n)(\tau+x)^{-1}-a(\tau+x)^{\alpha-1}=0,
\end{equation}
which, for $\alpha\ne1/2,1/3$, involves finding roots of polynomials
of degree greater than 4. When $\alpha=1/2$ we have
\begin{eqnarray*}
x^{*}_{D}
&=& \frac{(k-2)^{2}}{12a^{2}}-\frac{\tau}{3}\\
&&{}+\frac{48a^{2}\tau
(n-1)(k-2)+[4a^{2}\tau-(k-2)^{2}]^{2}}{6\cdot
2^{1/3}a^{2}p_{1,D}(a,\tau
,n,k)}\\
&&{} +\frac{p_{1,D}(a,\tau,n,k)}{12\cdot2^{1/3}a^{2}},
\end{eqnarray*}
where
\begin{eqnarray*}
\hspace*{-4pt}&&p_{1,D}(a,\tau,n,k)\\
\hspace*{-4pt}&&\qquad=\bigl\{p_{2,D}+\bigl[p_{2,D}^{2}+4\bigl(-48a^{2}\tau(n-1)(k-2)-
\bigl(4a^{2}\tau-(k-2)^{2}\bigr)^{2}\bigr)^{3}\bigr]^{{1}/{2}}\bigr\}
^{{1}/{3}}
\end{eqnarray*}
and $p_{2,D}=p_{2,D}(a,\tau,n,k)$ with
\begin{eqnarray*}
p_{2,D}(a,\tau,n,k)
&=& 2(k-2)^{3}[(k-2)^{3}-12a^{2}\tau(k+4-6n)]\\
&&{} +96a^{4}\tau^{2}\bigl(k(k+2)+10-6n(k+4)+18n^{2}\bigr)\\
&&{}-128a^{6}\tau^{3}.
\end{eqnarray*}
Similarly, one finds that
\begin{eqnarray*}
x^{*}_{N}
&=&\frac{(k-2)^{2}}{12a^{2}}-\frac{\tau}{3}\\
&&{}+\frac{48a^{2}\tau
n(k-2)+[4a^{2}\tau-(k-2)^{2}]^{2}}{6\cdot2^{2/3}a^{2}p_{1,N}(a,\tau
,n,k)}\\
&&{}+\frac{p_{1,N}(a,\tau,n,k)}{12\cdot2^{1/3}a^{2}},
\end{eqnarray*}
where
\[
p_{1,N}(a,\tau,n,k)=
[p_{2,N}+(p_{3,N})^{{1}/{2}}]^{{1}/{3}}
\]
with
\[
p_{2,N}(a,\tau,n,k)=p_{2,D}(a,\tau,n,k)-144a^{2}\tau[4a^{2}\tau
(k+1-6n)-(k-2)^{3}]
\]
and
\begin{eqnarray*}
&&p_{3,N}(a,\tau,n,k)\\
&&\qquad= -3^{3}2^{10}a^{6}n^{2}\tau^{3}\bigl[(k-2)^{4}-2n(k-2)^{3}\\
&&\hspace*{67.5pt}\qquad\quad{} -4a^{2}\tau[8+2k^{2}+9n(3n+4)-2k(9n+4)]
+16a^{4}\tau^{2}\bigr].
\end{eqnarray*}
When $k\approx s n^{\alpha}$ [cf. (\ref{GG-alpha-diversity}) above] and
$\alpha=1/2$, it can be checked that
%
\begin{equation}\label{p1-convergence}
n^{-1}p_{1,i}(a,\tau,n,k)\rightarrow2^{1/3}s^{2},\qquad
i=N,D,
\end{equation}
from which
%
\begin{equation}\label{mode-convergence}
n^{-1}x^{*}_{i}\rightarrow(s/2a)^{2},\qquad
i=N,D.
\end{equation}
Using this fact, one finds that
\[
\frac{f_{N}(x^{*}_{N})}{f(x^{*}_{D})}
\approx\frac{1}{1+\tau(s/2a)^{-2}}.
\]
Computing also the ratio of the two second derivatives, it can be seen
that $w(n,k)\rightarrow1$,
which by means of (\ref{decompos-g0}) and (\ref{constraint})
implies
(\ref{pesi-stable-case}).
In order to find the speed at which $w(n,k)$ goes to 1, consider
\[
1-w(n,k)=\int_{0}^{\infty}\frac{\tau}{\tau+x}f(x)\,\dd x\Big/\int
_{0}^{\infty}f(x)\,\dd x.
\]
The denominator is unchanged, while the mode of $\tau/(\tau+x)f(x)$ is
\begin{eqnarray*}
x'^{*}_{N}
&=& \frac{(k-4)^{2}}{12a^{2}}-\frac{t}{3}\\
&&{}-\frac{-48a^{2}t(n-1)(k-4)-[4a^{2}t-(k-4)^{2}]^{2}}{6\cdot
2^{2/3}a^{2}p'_{1,D}(a,\tau,n,k)}\\
&&{} +\frac{p'_{1,D}(a,\tau,n,k)}{12\cdot2^{1/3}a^{2}},
\end{eqnarray*}
where
\begin{eqnarray*}
\hspace*{-3pt}&&p'_{1,D}(a,\tau,n,k)\\
\hspace*{-3pt}&&\qquad=\bigl\{p_{2,D}+\bigl[p_{2,D}^{2}+4\bigl(-48a^{2}t(n-1)(k-4)-
\bigl(4a^{2}t-(k-4)^{2}\bigr)^{2}\bigr)^{3}\bigr]^{{1}/{2}}\bigr\}^{{1}/{3}}
\end{eqnarray*}
and $p_{2,D}=p_{2,D}(a,\tau,n,k)$ with
\begin{eqnarray*}
&&p'_{2,D}(a,\tau,n,k)\\
&&\qquad= 2(k-4)^{3}[(k-4)^{3}-12a^{2}\tau(k+2-6n)]\\
&&\qquad\quad{} +96a^{4}\tau^{2}\bigl(k(k-2)+10-6n(k+2)+18n^{2}\bigr)
-64a^{6}\tau^{3}.
\end{eqnarray*}
Moreover, $p'_{2,D}$ satisfies (\ref{p1-convergence}), and (\ref
{mode-convergence}) follows. Unfortunately the fact that the two modes
grow with an\vadjust{\goodbreak} asymptotically equivalent rate is too rough an
approximation for our purposes here, which ignores how far apart they
are if this is negligible with respect to the growth speed. Indeed, it
turns out that
\[
n^{-1/2}(x^{*}_{D}-x^{*}_{N})\rightarrow s/a^{2}.
\]
Using this information in the Laplace approximation for $w(n,k)-1$ yields
\[
n\bigl(1-w(n,k)\bigr)\rightarrow2a\sqrt{\tau}/s=\beta/s
\]
with $\beta$ as in the statement of the proposition.
From (\ref{decompos-g0}) it is now easy to see that
\[
ng_{0}(n,k)-\alpha k w(n,k)=n\bigl(1-w(n,k)\bigr),
\]
which provides the second order approximation for $g_{0}(n,k)$,
\[
g_{0}(n,k)=\frac{\alpha
k}{n}+\frac{\beta/s_{n}}{n}+o(n^{-1}),
\]
where the first term is of order $n^{-1/2}$, yielding immediately, by
means of (\ref{constraint}), the second order term for $g_{1}(n,k)$,
that is,
\[
g_{1}(n,k)=\frac{1}{n}-\frac{\beta/s_{n}}{n^{2}}+o(n^{-2}).
\]
%

\subsection*{Proof of Theorem \protect\ref
{theoremparticle-convergence}}
We can write the generator of (\ref{[K,Z]}) as
%
\begin{equation}\label{A-n}
\mathcal{B}_{n}(f_{0}\times f_{1})=f_{1}\mathcal
{B}_{n,0}f_{0}+f_{0}\mathcal{B}_{n,1}f_{1},
\end{equation}
where $\mathcal{B}_{n,0}$ and $\mathcal{B}_{n,1}$ drive $K_{n}(\cdot
)/n^{\alpha}$ and $w(X^{(n)}(\cdot))$, respectively, and $f_{0}\in
\mathscr{D}(\mathcal{A}_{0})$,
with $\mathscr{D}(\mathcal{A}_{0})$ as in (\ref{domain-A0}), $
f_{1}\in C^{2}(\nabla_{n})$ and $\nabla_{n}$ is as in (\ref{nabla-n}).
%
Based on~(\ref{K-transition}), while retaining $m_{1,n}$ temporarily, we can write
\begin{eqnarray*}
&&\mathcal{B}_{n,0}f_{0}\biggl(\frac{k}{n^{\alpha}}\biggr)\\
&&\qquad= \lim_{\delta t\downarrow0}\frac{1}{\delta t}
\biggl\{\delta t\biggl[
f_{0}\biggl(\frac{k+1}{n^{\alpha}}\biggr)
\biggl(1-\frac{m_{1,n}}{n}\biggr)g_{0}^{(n-1,k)}\\
&&\qquad\quad\hspace*{48pt}{} +f_{0}\biggl(\frac{k-1}{n^{\alpha}}\biggr)
\frac{m_{1,n}}{n}g_{1}^{(n-1,k-1)}\bigl(n-1-\alpha(k-1)\bigr) \\
&&\qquad\quad\hspace*{48pt}{} +f_{0}\biggl(\frac{k}{n^{\alpha}}\biggr)\biggl[1-\biggl(1-\frac
{m_{1,n}}{n}\biggr)g_{0}^{(n-1,k)}\\
&&\qquad\quad\hspace*{102pt}{} - \frac{m_{1,n}}{n}g_{1}^{(n-1,k-1)}\bigl(n-1-\alpha(k-1)\bigr)\biggr]\biggr]\\
&&\qquad\quad\hspace*{87.7pt}{} +(1- \delta t)f_{0}\biggl(\frac{k}{n^{\alpha}}\biggr)+O((\delta
t)^{2})-f_{0}\biggl(\frac{k}{n^{\alpha}}\biggr)\biggr\}\\
&&\qquad=
\biggl(1-\frac{m_{1,n}}{n}\biggr)g_{0}^{(n-1,k)}
\biggl[f_{0}\biggl(\frac{k+1}{n^{\alpha}}\biggr)-f_{0}\biggl(\frac
{k}{n^{\alpha}}\biggr)\biggr] \\[-2pt]
&&\qquad\quad{} + \frac{m_{1,n}}{n}g_{1}^{(n-1,k-1)}\bigl(n-1-\alpha(k-1)\bigr)
\biggl[f_{0}\biggl(\frac{k-1}{n^{\alpha}}\biggr)-f_{0}\biggl(\frac
{k}{n^{\alpha}}\biggr)\biggr].
\end{eqnarray*}
An application of Taylor's theorem yields
\begin{eqnarray*}
&&\mathcal{B}_{n,0}f_{0}\biggl(\frac{k}{n^{\alpha}}\biggr)\\[-2pt]
&&\qquad=
\biggl(1-\frac{m_{1,n}}{n}\biggr)g_{0}^{(n-1,k)}
\\[-2pt]
&&\qquad\quad{}\times
\biggl[f_{0}\biggl(\frac{k}{n^{\alpha}}\biggr)+\frac{1}{n^{\alpha
}}f_{0}'\biggl(\frac{k}{n^{\alpha}}\biggr)\\[-2pt]
&&\qquad\quad\hspace*{17pt}{}+\frac{1}{2n^{2\alpha
}}f_{0}''\biggl(\frac{k}{n^{\alpha}}\biggr)+O(n^{-3\alpha})
-f_{0}\biggl(\frac{k}{n^{\alpha}}\biggr)\biggr] \\[-2pt]
&&\qquad\quad{} + \frac{m_{1,n}}{n}g_{1}^{(n-1,k-1)}\bigl(n-1-\alpha(k-1)\bigr)
\\[-2pt]
&&\hspace*{11pt}\qquad\quad{}\times
\biggl[f_{0}\biggl(\frac{k}{n^{\alpha}}\biggr)-\frac{1}{n^{\alpha
}}f_{0}'\biggl(\frac{k}{n^{\alpha}}\biggr)+\frac{1}{2n^{2\alpha
}}f_{0}''\biggl(\frac{k}{n^{\alpha}}\biggr)\\[-2pt]
&&\qquad\quad\hspace*{106pt}{}+O(n^{-3\alpha})
-f_{0}\biggl(\frac{k}{n^{\alpha}}\biggr)\biggr] \\[-2pt]
&&\qquad=
\frac{1}{n^{\alpha}}f_{0}'\biggl(\frac{k}{n^{\alpha}}\biggr)
\biggl[\biggl(1-\frac{m_{1,n}}{n}\biggr)g_{0}^{(n-1,k)}
\\[-2pt]
&&\qquad\quad\hspace*{56pt}{}
-\frac{m_{1,n}}{n}g_{1}^{(n-1,k-1)}\bigl(n-1-\alpha(k-1)\bigr)\biggr] \\[-2pt]
&&\qquad\quad{}
+\frac{1}{2n^{2\alpha}}f_{0}''\biggl(\frac{k}{n^{\alpha}}\biggr)
\biggl[\biggl(1-\frac{m_{1,n}}{n}\biggr)g_{0}^{(n-1,k)}
\\[-2pt]
&&\qquad\quad\hspace*{82pt}{}
+\frac{m_{1,n}}{n}g_{1}^{(n-1,k-1)}\bigl(n-1-\alpha(k-1)\bigr)
\biggr]\\[-2pt]
&&\qquad\quad{}+O(n^{-4\alpha}).
\end{eqnarray*}
From (\ref{exp-increm}) and (\ref{exp-variance}) we have
\begin{eqnarray*}
\mathcal{B}_{n,0}f_{0}\biggl(\frac{k}{n^{\alpha}}\biggr)
&=&
\frac{1}{n^{\alpha}}f_{0}'\biggl(\frac{k}{n^{\alpha}}\biggr)
\biggl[\frac{\beta/s_{n}}{n}+o(n^{-1})\biggr] \\[-2pt]
&&{} +\frac{1}{2n^{2\alpha}}f_{0}''\biggl(\frac{k}{n^{\alpha}}\biggr)
\biggl[\frac{s_{n}}{n^{\alpha}}+O(n^{-1-2\alpha})\biggr]+O(n^{-4\alpha})\\[-2pt]
&=&
\frac{\beta/s_{n}}{n^{1+\alpha}}f_{0}'\biggl(\frac{k}{n^{\alpha}}\biggr)
+\frac{s_{n}}{2n^{3\alpha}}f_{0}''\biggl(\frac{k}{n^{\alpha}}
\biggr)+o(n^{-1-\alpha}),
\end{eqnarray*}
from which it follows that
%
\begin{equation}\label{A0-A0->0}
\|\mathcal{A}_{0}f_{0}-n^{3/2}\mathcal{B}_{n,0}f_{0}\|\rightarrow0
\end{equation}
as $n\rightarrow\infty$, with $\mathcal{A}_{0}$ as in (\ref{S-generator}).
Since (\ref{conv-distr-iniz}) implies $K_{n}(0)/n^{\alpha}\Rightarrow
S_{0}$, the previous expression, together with Theorems 1.6.1 and 4.2.5
in \citet{EK86}, implies (\ref{weak-conv-S}) with $C_{[0,\infty
)}([0,\infty))$ replaced by $D_{[0,\infty)}([0,\infty))$, while~(\ref
{weak-conv-S}) follows from relativization of the Skorohod topology
to\break
$C_{[0,\infty)}([0,\infty))$.

In order to describe $\mathcal{B}_{n,1}$ in (\ref{A-n}), define
%
\begin{equation}\label{phi}
\phi_{n}(\mu)=f_{1}(\langle h_{1},\mu\rangle,\ldots,\langle
h_{n},\mu\rangle),\qquad
f_{1}\in C_{0}^{2}(\mathbb{R}^{n}), h_{i}\in C(\mathbb{X}),
\end{equation}
and
%
\begin{equation}\label{empirical-process}
\mu_{n}(t)=\frac{1}{n}\sum_{i=1}^{n}\delta_{X_{i}(t)},\qquad
t\ge0.
\end{equation}
Then the generator of the $\mathscr{P}(\mathbb{X})$-valued process
$\mu
_{n}(\cdot)=\{\mu_{n}(t),t\ge0\}$ can be written
%
\begin{eqnarray}\label{m-v-gen-1}
\mathbb{B}_{n}\phi_{n}(\mu)
&=& \frac{1}{n}\sum_{i=1}^{n}g_{0}^{(n-1,i)}
[\langle Ph_{i},\mu\rangle-\langle h_{i},\mu\rangle] \,\frac
{\partial
f_{1}}{\partial z_{i}}\nonumber\\
&&{}- \frac{\alpha}{n}\sum_{i=1}^{n}g_{1}^{(n-1,i)}K_{n-1,i}
\bigl[\bigl\langle Q^{(n-1,i)}h_{i},\mu\bigr\rangle-\langle h_{i},\mu\rangle\bigr]
\,\frac
{\partial f_{1}}{\partial z_{i}}\\
&&{}+\frac{1}{n}\sum_{1\le k\ne i\le n}g_{1}^{(n-1,i)}[\langle
h_{i}h_{j},\mu\rangle-\langle h_{i},\mu\rangle\langle h_{j},\mu
\rangle] \,\frac{\partial^{2} f_{1}}{\partial z_{i}\,\partial
z_{j}},\nonumber
\end{eqnarray}
where $g_{z_{i}}$ is the derivative of $g$ with respect to its $i$th
argument and $Q^{(n-1,i)}$ is defined as
%
\begin{equation}\label{Q}
Q^{(n-1,i)}g(z)=\int g(y)p_{n-1,i}^{*}(z,\dd y),\qquad
g\in C_{0}(\mathbb{X}),
\end{equation}
with
%
\begin{equation}\label{p*}
p_{n-1,i}^{*}(z,\dd y)=\frac{1}{K_{n-1,i}}\sum
_{j=1}^{K_{n-1,i}}\delta
_{x_{j}^{*}}(\dd y).
\end{equation}
Unlabel now the model by choosing $\phi_{n}(\mu)$ as in (\ref{phi})
with $h_{j}(\cdot)$ being the indicator function of the $j$th largest
atom in $\mu$, so that $\langle h_{j},\mu\rangle=z_{j}$ is the relative
frequency associated with the $j$th most abundant species. Note that
some arguments of $f_{1}(z_{1},\ldots,z_{n})$ can be null since
$K_{n-1,i}\le K_{n}\le n$.
With this choice we have\vadjust{\goodbreak} $\langle h_{i}h_{j},\mu\rangle-\langle
h_{i},\mu\rangle\langle
h_{j},\mu\rangle=z_{i}\delta_{ij}-z_{i}z_{j}$, where $\delta_{ij}$
is the
Knonecker delta, and, under (\ref{mutation}) and (\ref{Q}),
\begin{eqnarray*}
\langle Ph_{i},\mu\rangle-\langle h_{i},\mu\rangle
&=&\sum_{j=1}^{n}p_{1}(x_{j}^{*},\dd x_{i}^{*})z_{j}-z_{i}\\[-2pt]
&=&\sum_{j=1}^{n}[p_{1}(x_{j}^{*},\dd x_{i}^{*})-\delta_{ij}]z_{j}
\end{eqnarray*}
and
\[
\bigl\langle Q^{(n-1,i)}h_{i},\mu\bigr\rangle-\langle h_{i},\mu\rangle
=\sum_{j=1}^{n}[p_{n-1,i}^{*}(x_{j}^{*},\dd x_{i}^{*})-\delta
_{ij}]z_{j}.
\]
It follows that $\mathbb{B}_{n}$ reduces to
%
\begin{eqnarray}\label{m-v-gen-2}
\mathcal{B}_{n,1}f_{1}
&=&\frac{1}{n}\sum_{i=1}^{n}g_{0}^{(n-1,i)}
\Biggl(\sum_{j=1}^{n}[p_{1}(x_{j}^{*},\dd x_{i}^{*})-\delta_{ij}
]z_{j}\Biggr)\,\frac{\partial f_{1}}{\partial z_{i}}\nonumber\\[-2pt]
&&{} -\frac{\alpha}{n} \sum_{i=1}^{n}g_{1}^{(n-1,i)}K_{n-1,i}\Biggl(\sum
_{j=1}^{n}[p_{n-1,i}^{*}(x_{j}^{*},\dd x_{i}^{*})-\delta_{ij}
]z_{j}\Biggr)\,\frac{\partial f_{1}}{\partial z_{i}}\\[-2pt]
&&{} +\frac{1}{n}\sum_{i,j=1}^{n}g_{1}^{(n-1,i)}z_{i}(\delta
_{ij}-z_{j})\,\frac{\partial^{2}f_{1}}{\partial z_{i}\,\partial
z_{j}}.\nonumber
\end{eqnarray}
Here a mutation from type $i$ to type $j$ occurs at rate
$q_{ij}=q_{ij}(z)$ given by
%
\begin{eqnarray}\label{mutationrate}
q_{ij}&=&\frac{1}{n}\bigl[g_{0}^{(n-1,i)}
[p_{1}(i,\{j\})-\delta_{ij}]\nonumber\\[-8pt]\\[-8pt]
&&\hspace*{10pt}{}-\alpha g_{1}^{(n-1,i)}K_{n-1,i}[p_{n-1,i}^{*}(i,\{j\})-\delta
_{ij}]\bigr],\nonumber
\end{eqnarray}
where $p_{1}(i,\{j\})$ and $p_{n-1,i}^{*}(i,\{j\})$ stand for
$p_{1}(x_{i}^{*},\dd x_{j}^{*})$ and $p_{n-1,i}^{*}(x_{i}^{*},\dd x_{j}^{*})$.
When (\ref{nu0}) holds, from the nonatomicity of $\nu_{0}$ we have
$p_{1}(z,\dd y)=0$ for every $y\in\mathbb{X}$, and when (\ref{p*})
holds, we have $p_{n-1,i}^{*}(z,\dd y)=K_{n-1,i}^{-1}$, from which
%
\begin{eqnarray}\label{mutation-rate2}
q_{ij} &=& \frac{1}{n}
\bigl[-\delta_{ij}g_{0}^{(n-1,i)}-\alpha g_{1}^{(n-1,i)}+\alpha
g_{1}^{(n-1,i)}K_{n-1,i}\delta_{ij}\bigr]\nonumber\\[-2pt]
&=& \frac{1}{n}\bigl[-\delta_{ij}\bigl(1-(n-1)g_{1}^{(n-1,i)}\bigr)-\alpha
g_{1}^{(n-1,i)}\bigr]\nonumber\\[-2pt]
&=& \frac{1}{n}\biggl[-\delta_{ij}\biggl(1-(n-1)\biggl(\frac{1}{n-1}-\frac
{\beta/s}{(n-1)^{2}}\biggr)\biggr)\\[-2pt]
&&\hspace*{79.6pt}{} -\alpha\biggl(\frac{1}{n-1}-\frac{\beta/s}{(n-1)^{2}}\biggr)\biggr]
+o(n^{-2})\nonumber\\[-2pt]
&=& -\frac{\delta_{ij}\beta/s}{n(n-1)}-\frac{\alpha
}{n(n-1)}+o(n^{-2}),\nonumber
\end{eqnarray}
where the second equality follows from (\ref{constraint}) and the third
from (\ref{g1-approximation2}). Once again it is clear that the key
point for determining the limiting behavior of the diffusion is the
second order approximation of the predictive weights, as obtained in
Proposition \ref{second-order-approximation}.
Hence, we have
\[
\sum_{i=1}^{n}\sum_{j=1}^{n}q_{ij}z_{j}
=-\sum_{i=1}^{n}\biggl[\frac{z_{i}\beta/s}{n(n-1)}+\frac{\alpha
}{n(n-1)}+o(n^{-2})\biggr],
\]
from which (\ref{m-v-gen-2}), substituting (\ref{g1-approximation2}) in
the third term, reduces to
\begin{eqnarray*}
\mathcal{B}_{n,1}f_{1}
&=& \sum_{i,j=1}^{n}[n^{-2}-O(n^{-3})]z_{i}(\delta
_{ij}-z_{j})\,\frac{\partial^{2}f_{1}}{\partial z_{i}\,\partial z_{j}}\\
&&{} -\sum_{i=1}^{n}\biggl[\frac{z_{i}\beta/s}{n(n-1)}+\frac{\alpha
}{n(n-1)}+o(n^{-2})\biggr]
\,\frac{\partial f_{1}}{\partial z_{i}},
\end{eqnarray*}
%
%
%
which in turn implies that
%
\begin{equation}\label{A1-A1->0}
\|\mathcal{A}_{1}f_{1}-(n^{2}/2)\mathcal{B}_{n,1}\tilde
r_{n}f_{1}\|\rightarrow0
\end{equation}
with $\mathcal{A}_{1}$ as in (\ref{gen-decomposition}), $f_{1}\in
\mathscr{D}(\mathcal{A}_{1})$ as in (\ref{domain-A-1}) and $\tilde
r_{n}$ as in (\ref{rn-tilde}).
From (\ref{A0-A0->0}) and (\ref{A1-A1->0}) it follows that
\[
\|\mathcal{A}(f_{0}\times f_{1})-(\tilde r_{n}f_{1})n^{3/2}\mathcal
{B}_{n,0}f_{0}-f_{0}n^{2}\mathcal{B}_{n,1}\tilde
r_{n}f_{1}/2\|\rightarrow0
\]
with $\mathcal{A}$ as in (\ref{operatorgibbs}) and $(f_{0}\times
f_{1})$ as in (\ref{domain-B}). The fact that (\ref
{eqparticle-convergence}) holds in $D_{[0,\infty)\times\overline
\nabla
_{\infty}}([0,\infty))$ follows from\vspace*{1pt} the density of (\ref{domain-B}) in
$C_{0}([0,\infty)\times\overline\nabla_{\infty})$, together with
Theorems 1.6.1 and 4.2.11 in \citet{EK86}, which imply, respectively,
the strong convergence of the associated semigroups on $C_{0}([0,\infty
)\times\overline\nabla_{\infty})$, similarly to (\ref
{eqSG-convergence}), and the weak convergence of the probability
measures induced on $D_{[0,\infty)\times\overline\nabla_{\infty
}}([0,\infty))$. The same assertion with $C_{[0,\infty)\times
\overline
\nabla_{\infty}}([0,\infty))$ in place of $D_{[0,\infty)\times
\overline
\nabla_{\infty}}([0,\infty))$ now follows from relativization of the
Skorohod topology.
\end{appendix}

\section*{Acknowledgments}

The authors would like to thank an Associate Editor and a Referee for
their extremely careful reading of previous versions of the paper and
for providing constructive suggestions that greatly improved the paper,
and the Editor for his support.
%
%
Special thanks also go to Anshui Li, Antonio Lijoi, Bertrand Lods and Igor
Pr\"unster.



\printaddresses

\end{document}